\documentclass{amsart}
\makeindex  

\usepackage{amsmath,amssymb}
\usepackage{amsthm}
\usepackage{amscd}

\newcommand{\bQ}{{\mathbb Q}}

\newcommand{\bZ}{{\mathbb Z}}

\newcommand{\Gal}{{\rm Gal}}

\def\Cbb{{ C}}

\def\Nbb{{ N}}
\def\Pbb{{P}}
\def\Qbb{{ Q}}

\def\Zbb{{ Z}}

\def\Ic{{\mathcal I}}

\def\Kc{{\mathcal K}}
\def\Lc{{\mathcal L}}
\def\Mc{{\mathcal M}}
\def\Nc{{\mathcal N}}

\def\zfk{{\mathfrak z}}
\def\pfk{{\frak p}}

\def\01{{\overrightarrow{01}}}
\def\10{{\overrightarrow{10}}}

\def\hpb{\hfill $\Box$}

\def\al{\alpha}
\def\be{\beta}
\def\de{\delta}
\def\ga{\gamma}
\def\si{\sigma}

\def\lam{\lambda}
\def\Ga{\Gamma}

\def\Aut{{\rm Aut}}

\def\log{{\rm log}}

\def\pd1{{\partial \Delta [1]}} 
\def\d1{{\Delta [1]}}

\def\sib{{\langle \sigma \rangle }}
\def\sibn{{\langle \sigma \rangle _n }}

\def\xn{\xi _N}
\def\nn{\frac{1}{N}}

\def\Cbb{{\mathbb C}}

\def\Nbb{{\mathbb N}}
\def\Pbb{{\mathbb P}}
\def\Qbb{{\mathbb Q}}

\def\Zbb{{\mathbb Z}}

\def\Hc{{\mathcal H}}
\def\Ic{{\mathcal I}}

\def\Lc{{\mathcal L}}
\def\Mc{{\mathcal M}}
\def\Nc{{\mathcal N}}

\def\Yc{{\mathcal Y}}
\def\Zc{{\mathcal Z}}

\def\z2{{{\Zbb [{\frac{1}{2}}]}}}

\def\ffk{{\mathfrak f}}
\def\pfk{{\mathfrak p}}

\def\tf{{\mathfrak t}}

\def\uf{{\mathfrak u}}
\def\sf{{\mathfrak s}}
\def\zf{{\mathfrak z}}

\def\al{{\alpha}}
\def\be{{\beta}}
\def\ga{{\gamma}}
\def\si{{\sigma}}
\def\de{{\delta}}

\def\sip{{ \sum _{i=0}^{p^n-1}}}
\def\sbp{{ \sum _{b=0}^{p^n-1}}}
\def\sap{{ \sum _{a=0}^{p^n-1}}}
\def\sbp{{ \sum _{b=0}^{p^n-1}}}
\def\sa{\sum _{a,b}}

\def\chs{{\chi (\si )}}
\def\ls{\langle \si \rangle _n}
\def\cs{\langle c \rangle _n}

\def\sb{{\langle \si \rangle }}
\def\sbn{{\langle \si \rangle _n }}

\def\zpn{\Zbb /(p^n)}
\def\zp{\Zbb _p}

\def\01{{\overset{\to}{01}}}
\def\10{\overset{\to}{10}}

\def\hpb{\hfill $\Box$}

\def\Aut{{\rm Aut}}
\def\jn{{\frac{1}{N}}}

\def\zi{{ \zp [[(\zp )^i]]}}

\def\bQ{{\overline \Qbb}}
\def\zp{{\Zbb _p}}
\def\qp{{\Qbb _p}}

\title{Octagonal relations}
\author{Zdzis{\l}aw Wojtkowiak}
\begin{document}

\date{\today}

\maketitle

\tableofcontents

\begin{abstract}
The absolute  Galois group $\Gal (\overline \Qbb /\Qbb )$ acts  on the  \'etale fundamental group of $\Pbb ^1_ {\overline \Qbb }\setminus (\{ 0,\infty \}\cup \mu _{{p^n}})$  and on the canonical path $\pi _n$ from $\01$ to $\frac{1}{p^n}\10$ on $\Pbb ^1_ {\overline \Qbb }\setminus (\{ 0,\infty \}\cup \mu _{{p^n}})$ for all $n\in \Nbb$.  We show the octagonal relation for the path  $\pi _n$  in the Galois context. The above Galois actions 
 leads to the natural construction of measures on $(\zp )^m$ for $m\in \Nbb$. We study the symmetries of such measures derived from the octagonal relations. We calculate explicitely the symmetry for $m=2$. We show also that the Galois action (pro-$p$) on the path $\pi _0$ from $\01$ to $\10$ on 
 $\Pbb ^1_ {\overline \Qbb }\setminus (\{ 0,1,\infty \} )$ is determined and determinds the measures mentioned above for all $n\in \Nbb$.
\smallskip

\noindent
R\'esum\'e.
Le groupe de Galois absolu $\Gal (\overline \Qbb /\Qbb )$ op\`ere sur le groupe fondamental \'etale de $\Pbb ^1_ {\overline \Qbb }\setminus (\{ 0,\infty \}\cup \mu _{p^n})$ et sur le canonique chemin 
 $\pi _n$ de  $\01$ \`a $\frac{1}{p^n}\10$ sur $\Pbb ^1_ {\overline \Qbb }\setminus (\{ 0,\infty \}\cup \mu _{{p^n}})$ pour tous $n\in \Nbb$.
On montre la relation octogonale pour le chemin  $\pi _n$  en Galois contexte. 
 Ces actions conduisent \` a des constructions  des mesures sur $(\zp )^m$ pour tous $m\in \Nbb$.
On \'etudie les sym\'etries de ces mesures d\'eduites des relations octogonales. On calcule explicitement la sym\'etrie pour $m=2$.
On montre aussi que l'action de Galois  (pro-$p$) sur le chemin $\pi _0$ de $\01$ to $\10$ sur
 $\Pbb ^1_ {\overline \Qbb }\setminus (\{ 0,1,\infty \} )$ est d\'etermin\'ee et d\'etermine les mesures mentionn\'ees au-dessus.
\end{abstract}

\section{Introduction}
\smallskip

The research of representations of Galois groups on non-abelian fundamental groups were started by A.  Grothendieck, Y. Ihara (see \cite{I}), P. Deligne (see \cite{D}). In this paper we study some aspects of measures associated to Galois actions  on  a family  of cyclic coverings of $\Pbb ^1\setminus \{ 0,1,\infty \}$.  The present  paper is a continuation of \cite{W10}.

We fix once for all an embedding
$
\overline \Qbb \subset \Cbb\;.
$
Let $N$ be a positive integer. We set
$
\mu _N:=\{ z\in \overline \Qbb \mid z^N=1\}.
$

Let $V:=\Pbb ^1_ {\overline \Qbb }\setminus (\{ 0,\infty \}\cup \mu _N)$ and let $\pi  : \01 \leadsto \frac{1}{N} \10$ be the canonical path from $\01$ to $ \frac{1}{N} \10$ along the real interval the interval $[0,1]$. 
Let $\si \in \Gal (\overline \Qbb /\Qbb )$. We set
$$
\ffk _\pi (\si ):=\pi ^{-1}\cdot \si (\pi )=\pi ^{-1}\cdot \si \cdot \pi \cdot \si ^{-1}\in  \pi _1 ^{\rm et}(V, \01 ).
$$
Below we shall describe briefly the main results of the paper. In section 2 we show that the element $ \ffk _\pi (\si )$ satisfies the octagonal relation
(see \cite[Proposition 5.3.]{W15}, \cite{E}, \cite[the equation (2,2), Theorem 3.1,]{EF}, \cite{G} for the complex analogues).

Let us fix a rational prime $p$.  In the rest of the paper we shall restrict our study to the Galois action on the pro-$p$ completion of \'etale fundamental groups and on pro-$p$ paths. As in other our papers we shall embed pro-$p$ groups in $\Qbb _p$-algebras of formal power series in non-commuting variables.
Starting from section 3, we change slightly the notation as we shall study the tower of coverings
 $
 \{ \Pbb ^1_ {\overline \Qbb }\setminus (\{ 0,\infty \}\cup \mu _{p^n})\} _{n\geq 0}
$
of $\Pbb ^1_ {\overline \Qbb }\setminus \{ 0,1,\infty  \}$.

Let $V_n:= \Pbb ^1_ {\overline \Qbb }\setminus (\{ 0,\infty \}\cup \mu _{p^n})$ for $n\geq 0$ and let 
$\pi _n : \01 \leadsto \frac{1}{p^n} \10$ be the canonical path from $\01$ to $ \frac{1}{p^n} \10$ along the real   interval $[0,1]$. 

In section 3 we recall from \cite{W11} the construction of measures $\Kc _r(\si )$ on $(\zp )^r$ for $\si \in  \Gal (\overline \Qbb /\Qbb )$ and $r=1,\, 2,\, 3 \ldots$ associated with the path $\pi _0$ and the tower of coverings $V_n$ of $V_0=\Pbb ^1_ {\overline \Qbb }\setminus \{ 0,1,\infty  \}$.
The measure $\Kc _1(\si )$ is (the pro-$p$ part of) the measure $D\log \Ga _\si \,  - \,(D\log \Ga _ \si )(0)$ studied by Anderson and Ihara (see \cite[Theorem 5 ]{A}, \cite[Theorem 2]{I1}, \cite[page 1]{I2}) and where $\Ga _\si$ is the Anderson hyperadelic gamma function.

The octagonal relations satisfied by the elements $\ffk _{\pi _n}(\si ):=\pi _n^{-1}\cdot \si (\pi _n )=\pi _n ^{-1}\cdot \si \cdot \pi  _n \cdot \si ^{-1}\in  \pi _1 ^{\rm et}(V_n, \01 )$ ($n=1,\,2,\ldots$) imply certain relations for measures $\Kc _r(\si )$ as these measures are formed from the coefficients of the elements  $\ffk _{\pi _n}(\si )$.
Our aim is to write down explicitely such relations, at least for small $r$.

In section 4 we shall write (almost) explicitely relations between coefficients in degrees 1 and 2, which are consequence of the octagone relation ( see Proposition 4.3). When trying to adjust these coefficients to get measures, some new measures on $\zp$ and $(\zp )^2$ appear naturally. These measures are presented in section 5. In section 6 we present some elementary operations on measures which also will be needed later.

In \cite[Theorem 2.6.]{W11} we have shown that the coefficients of the image of the element $ \ffk _{\pi _0}(\si )$ (and more generally  $\ffk _{\ga }(\si )$ for any path $\ga$) in $\Qbb _p\{ \{ X,Y\} \} $ can be expressed by integrals of polynomials againt measures $\Kc _r(\si )$.
In section 7 we shall show that   the coefficients of the image of  $\ffk _{\pi _0}(\si )$ in  $\Qbb _p\{ \{ X,Y\} \} $determine uniquely measures  $\Kc _r(\si )$
for $r=1,\, 2, \ldots  \,$ .  We shall also indicate that these two results are purely group theoretical.

In section 8 we write down the relations satisfied by the measures  $\Kc _1(\si )$ and   $\Kc _2(\si )$ and which are the consequence of the octagonal  relations (see Proposition  8.1 and Theorems 8.2 and 8.4).

Let $m>1$. We were hoping to recover from measures $\Kc _m(\si )$   $p$-adic multiple zeta functions. We have seen in \cite{F}, that the authors used the measure $E_{1,c}$ with $c\in \zp^\times$ to construct $p$-adic  multiple zeta functions. There is in fact  a way to relate the measure  $\Kc _m(\si )$ to the measure  $E_{1, \chs }$. This result is presented in section 9.

In section 11 we state without proof analogues of some formulas from the present paper in the Betti - De Rham realisation. However the Betti - De Rham analogue of the formula from Theorem 8.4 is not presented. We shall try to write it in some future paper.

\bigskip

\section{Octagonal relation}
\smallskip

We recall that we fix once for all an embedding
$
\overline \Qbb \subset \Cbb
$
and that for a positive integer  $N$ 
$$
\mu _N:=\{ z\in \overline \Qbb \mid z^N=1\}\;\;{\rm and}\;\; \xi _N:=e^{\frac{2\pi i}{N}}.
$$
Let $V:=\Pbb ^1_ {\overline \Qbb }\setminus (\{ 0,\infty \}\cup \mu _N)$ and let  $V(\Cbb ):=\Pbb ^1 (\Cbb )\setminus (\{ 0,\infty \}\cup \mu _N)$.
Let $\pi  : \01 \leadsto \frac{1}{N} \10$ be the canonical path from $\01$ to $ \frac{1}{N} \10$ along the real interval  $[0,1]$.

For $a\in \Zbb /{(N)}$ let $R_a:V \to V$ be given by $R_a(\zfk ):=\xn ^a\zfk$. However because of certain asymmetry we shall take often indices between $0$ and $N$.

Let $s_a: {\overrightarrow{{0}{\xn ^a} } } \leadsto \01$ (resp. $t_a: {\overrightarrow{{ \infty}{1} } } \leadsto  {\overrightarrow{{\infty } {\xn ^a} }}$,
 resp. $u_a:\nn  {\overrightarrow { \xn ^a \infty  } } \leadsto \nn  {\overrightarrow{ \xn ^a 0  } }$)
be a path in an infinitesimal neighbourhood of $0$ (resp. $\infty$, resp. $\xn ^a$) --  an arc of a circle in the clockwise direction.

We choose generators of $\pi _1(V(\Cbb ),\01 )$ in the following way.
The element $x$ is a loop around $0$ in the positive direction  in an infinitesimal neighbourhood of $0$. For $0\leq a<N$  let $y^\prime _a \in \pi _1(V(\Cbb ),                       
 {\overrightarrow{{\xn ^a}{0} } }    )$ be a loop around $\xi _N$ in the positive direction in an  infinitesimal neighbourhood of $\xn ^a$. We set $y_a:=s_a\cdot R_a(\pi )^{-1} \cdot y^\prime _a\cdot  R_a(\pi ) \cdot s_a ^{-1}$ (see Picture 1).

\[
\,
\]

\[
\,
\]
$$
{\rm Picture\; 1}
$$

The elements $x,y_0,y_1,\ldots y_{N-1}$ are free generators of the free group $\pi _1(V(\Cbb ),\01 )$  and there are also free generators of the pro-finite free group $ \pi _1 ^{\rm et}(V, \01 )$.

 Let $k:V\to V$ be given by $k(\zfk )=1/\zfk$ and  let  $d:=k(\pi  )^{-1 }$. Let  $r:\nn\10\leadsto -\nn \10$ (=$ \nn {\overrightarrow {1\infty }}$)  be an arc in an infinitesimal neighbourhood of $1$ in the  clockwise direction from $\nn\10 $ to $-\nn\10$  ($=\nn {\overrightarrow {1\infty }}$). 
Let $\Ga :=d\cdot r\cdot \pi $.

Let $z_\infty \in \pi _1(V(\Cbb ),{\overrightarrow {\infty 1}} )$ be a loop around $\infty $ in an infinitesimal neighbourhood of $\infty$ in the positive direction.  Let us set $z:=\Ga ^{ -1}\cdot z_\infty \cdot \Ga $. Observe that 
\[
x\cdot  y_{N-1}\cdot \ldots   \cdot   y_{2}\cdot  y_{1}\cdot z\cdot  y_{0}=1
\]
in $\pi _1(V(\Cbb ),\01 )$.

For $a\in \Zbb /{(N)}$,  let  us set  $b_a:=R_a(\pi ^{-1}): \frac{1}{N} {\overrightarrow{{\xn ^a}{0 }} } \leadsto  {\overrightarrow{ 0{\xn ^a}} }$,
 $e:=k(\pi ):{\overrightarrow{ \infty 1} }\leadsto   \frac{1}{N} {\overrightarrow{  1 \infty } }$ and
 $c_a:=R_a(e): {\overrightarrow{{ \infty}{\xn ^a} } } \leadsto   \frac{1}{N}   {\overrightarrow{ {\xn ^a}{\infty  }} }$.

\medskip

The proofs of the next three lemmas we left to the reader as they are elementary though they required careful work.

\medskip

\noindent
{\bf Lemma 2.1.}
{\it  In the group  $ \pi _1 (V(\Cbb ), \01 )$  we have
\begin{enumerate}
\item[i)] $ \Ga ^{-1}\cdot k(x)\cdot \Ga =z$,

\item[ii)] $\Ga ^{-1}\cdot k(y_0)\cdot \Ga =y_0$,

\item[iii)] $\Ga ^{-1}\cdot k(y_k)\cdot \Ga =y_0\cdot x \cdot y_{N-1}\cdot \ldots y_{N-(k-1)}\cdot y_{N-k}\cdot (y_0\cdot x \cdot y_{N-1}\cdot \ldots y_{N-(k-1)} )^{-1}$ for $1\leq k \leq N-1$.
\end{enumerate}
}  \hpb

\medskip

Let us set $q_1:=r\cdot \pi$, $q_2:=d\cdot q_1$,  $q_3(a):=t_a\cdot q_2$, $q_4(a):=c_a\cdot q_3(a)$, $q_5(a):=u_a\cdot q_4(a)$, $q_6(a):=b_a\cdot q_5(a)$.

\medskip

\noindent
{\bf Lemma 2.2.}
{\it   In the group  $ \pi _1 (V(\Cbb ), \01 )$  we have
\begin{enumerate}
\item[i)] if     $0<a<N$  then $ q_3(a)^{-1}\cdot R_a(k(x))\cdot q_3(a) =z$,

\medskip

\item[ii)] if $0\leq b<a\leq N-1$ then $ q_3(a)^{-1}\cdot R_a(k(y_b))\cdot q_3(a) =(y_{a-b-1}\cdot \ldots \cdot y_2\cdot y_1)^{-1}\cdot y_{a-b} \cdot (y_{a-b-1}\cdot \ldots \cdot y_2\cdot y_1)  $,

\medskip

\item[iii)]  if $0 <a<b\leq  N-1$ then $ q_3(a)^{-1}\cdot R_a(k(y_b))\cdot q_3(a) =(y_{N+a-b-1}\cdot \ldots \cdot y_2\cdot y_1\cdot z)^{-1}\cdot y_{a-b} \cdot (y_{N+a-b-1}\cdot \ldots \cdot y_2\cdot y_1\cdot z)  $,

\medskip

\item[iv)]  $ q_3(a)^{-1}\cdot R_a(k(y_a))\cdot q_3(a) =y_0$.
\end{enumerate}
}  \hpb

\medskip

\noindent
{\bf Lemma 2.3.}
{\it   In the group  $ \pi _1 (V(\Cbb ), \01 )$  we have
\begin{enumerate}
\item[i)] $ q_6(a)^{-1}\cdot R_a(x)\cdot q_6(a)     =(y_{a-1}\cdot \ldots \cdot y_2\cdot y_1)^{-1}\cdot x \cdot (y_{a-1}\cdot \ldots \cdot y_2\cdot y_1)  $  for $0<a<N$,  

\medskip
  
\item[ii)]   $ q_6(a)^{-1}\cdot R_a(y_i)\cdot q_6(a) =(y_{a-1}\cdot \ldots \cdot y_2\cdot y_1)^{-1}\cdot y_{i+a} \cdot (y_{a-1}\cdot \ldots \cdot y_2\cdot y_1)  $ if $0\leq  a<N$, $0\leq i<N$ and $i+a<N$,

\medskip

\item[iii)]  $ q_6(a)^{-1}\cdot R_a(y_i)\cdot q_6(a) =(x\cdot y_{a-1}\cdot \ldots \cdot y_2\cdot y_1)^{-1}\cdot y_{i+a} \cdot (x\cdot y_{a-1}\cdot \ldots \cdot y_2\cdot y_1)  $    if  $0 <a<  N$, $0\leq i< N$ and $i+a\geq N$ 
\end{enumerate}
 with the convention that $y_{a-1}\cdot \ldots y_2\cdot y_1=1$ if $a=1$.}     \hpb

\bigskip

Let $a\in \hat \Zbb$. We define $\langle  a \rangle _N \in \Nbb$ by the following two conditions $\langle a \rangle _N \equiv a\; {\rm mod}\;  N$ and $0\leq \langle a \rangle _N <N$.
For $\si \in G_\Qbb$ we set 
\[
\sib =\sib _N:=\langle \chi (\si )\rangle _N.
\]

\medskip

Let $\si \in G_\Qbb$.    Then $\pi ^{-1}\cdot \si   ( \pi )=   \pi ^{-1}\cdot \si_{ \jn \10}\cdot \pi\cdot(\si_{\01})^{-1}  \in \pi _1 ^{\rm et}(V, \01 )$.     ( The path acts by analytic continuation and $\si _x$ denotes the action of $\si \in G_\Qbb$ on the fiber over $x$.)      Let us set
\[
\ffk _\pi(\si )=\ffk _\pi(\si )(x,y_0,\ldots ,y_{N-1}):= \pi ^{-1}\cdot \si_{ \jn \10}\cdot \pi\cdot(\si_{\01})^{-1}\; .
\]

\medskip

Observe that 
\[
\al :=s_1 \cdot b_1\cdot u_1\cdot c_1\cdot t_1\cdot d\cdot r\cdot \pi =1
\]
in    $\pi _1(V (\Cbb ),\01 )$, hence also  in   $\pi _1^{\rm et }(V,\01 )$ (see Picture 2). (See also \cite[Picture 4]{W10}, though the notation is different there.)

\[
\,
\]

\[
\,
\]

$$
{\rm Picture\; 2}
$$

\medskip

 Hence for any $\si \in G_\Qbb$,
\[
\al  ^{-1}\cdot \si _{\01  }\cdot \al   \cdot \si _{  \01} ^{-1   }=1 
\]
in  $\pi _1^{\rm et }(V,\01 )$.

\medskip

Let $\zf$ be the standard coordinate on $\Pbb _\bQ^1$. Then $\zf$ is a local parameter at $0$ corresponding to the tangential point $\01$.
Below we list local parameters and corresponding tangential points we need. (Notice that base points change covariantly and local parameters change contravariantly.)

For $0<k<N$, $\zf _k:=\xi _N^{-k}\zf$ corresponds to $\overrightarrow{0\xn ^k}$.
For $0\leq k <N$, $\uf _k := \xn ^k/\zf$ corresponds to $\overrightarrow{\infty \xn ^k}$, $\sf _k :=N\big( {\frac{\zf - \xn ^k}{\zf}}\big)$ corresponds to 
${\frac{1}{N}}\overrightarrow{\xn ^k \infty}$ and $\tf _k :=N \xn^{-k}(\xn ^k -\zf )$ corresponds to ${\frac{1}{N}}\overrightarrow{\xn ^k 0}$.

\medskip

\noindent
{\bf Lemma 2.4.}
{\it  Let $\si \in G_\Qbb$. Then in the group $\pi _1^{\rm et}(V,\01 )$  we have
\[
1=(s_1\cdot b_1\cdot u_1\cdot c_1\cdot t_1\cdot d\cdot r\cdot \pi )^{-1}\cdot s_{\sib }\cdot q_6(\sib )\cdot 
 q_6(\sib)^{-1 }\cdot \big( (s_{ \sib })^{-1}\cdot \si _{\01 }\cdot s_1\cdot  ( \si _ {\overrightarrow {0  {\xn ^1}  }   })^{ -1}\big) \cdot q_6(\sib )\cdot
\]
\[
q_5(\sib)^{-1 }\cdot \big( b_{\sib} ^{-1 }\cdot  \si _ {  \overrightarrow {0 {\xn ^1}  } }\cdot  b_1\cdot ( \si _ { \jn \overrightarrow {  {\xn ^1}   0  } })^{ -1}\big) \cdot  q_5(\sib) \cdot \,  q_4(\sib)^{-1 }\cdot \big( u_{\sib} ^{-1 }\cdot  \si _ { \jn  \overrightarrow { {\xn ^1}0      } }\cdot  u_1\cdot ( \si _ {  \jn \overrightarrow {  {\xn ^1}   \infty  } })^{ -1}\big)\cdot
 \]
\[
 q_4(\sib)\cdot    q_3(\sib)^{-1 }\cdot (c_{\sib} ^{-1 }\cdot  \si _ {  \jn  \overrightarrow { {\xn ^1}   \infty   } }\cdot  c_1\cdot ( \si _ {\overrightarrow { \infty {\xn ^1}     } })^{ -1})\cdot  q_3(\sib) \cdot
q_2^{-1}\cdot \big( t_{\sib} ^{-1 }\cdot  \si _ {\overrightarrow {\infty  {\xn ^1}      } }\cdot  t_1\cdot 
\]
\[
 ( \si _ {\overrightarrow { \infty  1     } })^{ -1})\cdot  q_2 \cdot
q_1^{-1}\cdot \big( d ^{-1 }\cdot  \si _ {\overrightarrow {\infty  1      } }\cdot  d\cdot ( \si _ { \jn  \overrightarrow {1 \infty     } })^{ -1}\big) \cdot  q_1 \cdot \pi ^{-1}\cdot \big( r^{-1}\cdot  \si _ {  \jn \overrightarrow {1 \infty     } }\cdot r \cdot ( \si _ {\jn  \overrightarrow {1 0     } })^{-1}\big) \cdot \pi \cdot
\]
\[
\big(\pi ^{-1}\cdot \si_{ \jn \10}\cdot \pi\cdot(\si_{\01})^{-1}\big) \; .
\]
}

\medskip

\noindent
{\bf  Proof.}   Canceling carefuly adjacent terms we see that the right hand side is $\al  ^{-1}\cdot \si _{\01  }\cdot \al   \cdot \si _{  \01} ^{-1   }=1$.
\hpb

\medskip

We shall express elements $q_1^{-1}\cdot \big( d ^{-1 }\cdot  \si _ {\overrightarrow {\infty  1      } }\cdot  d\cdot ( \si _ { \jn  \overrightarrow {1 \infty     } })^{ -1}\big) \cdot  q_1$,

\noindent
$q_3(\sib)^{-1 }\cdot (c_{\sib} ^{-1 }\cdot  \si _ {  \jn  \overrightarrow { {\xn ^1}   \infty   } }\cdot  c_1\cdot ( \si _ {\overrightarrow { \infty {\xn ^1}     } })^{ -1})\cdot  q_3(\sib) $ and 
$q_5(\sib)^{-1 }\cdot \big( b_{\sib} ^{-1 }\cdot  \si _ {  \overrightarrow {0 {\xn ^1}  } }\cdot  b_1\cdot ( \si _ { \jn \overrightarrow {  {\xn ^1}   0  } })^{ -1}\big) \cdot  q_5(\sib) $
by the element $\ffk _\pi(\si )$. 

\medskip

\noindent
In order to simplify the notation we set
\[
Q_1:=\jn \10  \, ,\; Q_2:= \jn \overrightarrow {1 \infty}  \; \; {\rm and}\; \;  Q_3(a):= \jn \overrightarrow {\xi _N^a \infty} \, ,\;  Q_4(a):= \jn \overrightarrow {\xi _N^a 0}\, 
\]
for $0<a<N$.

Observe that local parameters corresponding to tangential points $\jn \overrightarrow {1\infty }$ and $-\jn \10$ are the same.
\medskip

\noindent
{\bf Lemma 2.5.}
{\it Let $\si \in G_\Qbb$. Then in the group $\pi _1^{\rm et}(V,\01 )$ we have
\begin{enumerate}
\item[i)] $q_1^ {-1  }\cdot(d^{-1  }\cdot \si _ {\overrightarrow {   \infty 1  } }\cdot d\cdot (\si _ {Q_2})^{-1  })\cdot q_1= q_2^{-1  }\cdot k_*\big( \ffk _\pi (\si )\big)^{-1} \cdot q_2$,

\medskip

\item[ii)] $q_3(\sib)^{-1 }\cdot (c_{\sib} ^{-1 }\cdot  \si _ {Q_3(1) }\cdot  c_1\cdot ( \si _ {\overrightarrow { \infty {\xn ^1}     } })^{ -1})\cdot  q_3(\sib) =$
\newline
$q_3(\sib)^{-1 }\cdot \big( (R_{\sib})_*\big( k_*\big( \ffk _\pi (\si )\big)\big) \big) \cdot q_3(\sib) $,

\medskip

\item[iii)] $q_5(\sib)^{-1 }\cdot \big( b_{\sib} ^{-1 }\cdot  \si _ {\overrightarrow {0 {\xn ^1}  } }\cdot  b_1\cdot ( \si _ {Q_4(1) })^{ -1}\big) \cdot  q_5(\sib) =$
\newline
$q_6(\sib)^{-1 }\cdot \big( (R_{\sib})_*\big(  \ffk _\pi (\si )\big)^{-1  } \big) \cdot q_6(\sib) \, ,$

\medskip

\item[iv)]  $(s_1\cdot b_1\cdot u_1\cdot c_1\cdot t_1\cdot d\cdot r\cdot \pi )^{-1}\cdot s_{\sib }\cdot q_6(\sib )= y_   {\sib  -1}\cdot \ldots \cdot y_2\cdot y_1$.
\end{enumerate}
}

\smallskip

\noindent
{\bf  Proof.} We start with the proof of i). Observe that $k^*(\sf _0)=\tf _0$ and $k^*(\uf _0)=\zfk$.
We have the following commutative diagram of field isomorphisms

$$
\begin{matrix}

 \bQ ((\sf _0^{1/\infty})) &{\overset{(\si _ {Q_2})^{-1} }{\longrightarrow}} &     \bQ ((\sf _0^{1/\infty})) &{ \overset{k(\pi )^{-1}} {\longrightarrow} }&  \bQ ((\uf _0^{1/\infty})) &  { \overset{    \si    _{      \overrightarrow {\infty  1}    }                    } {\longrightarrow} } &    \bQ ((\uf _0^{1/\infty})) &   { \overset{k(\pi )} {\longrightarrow} }& \bQ ((\sf _0^{1/\infty}))  \\ \\
 {k^*  }\downarrow &&{k^* }\downarrow    && {k^* }\downarrow    && {k^* }\downarrow  &&  {k^* }\downarrow  \\ \\
 \bQ ((\tf _0^{1/\infty}))    &{\overset{ (\si _ {Q_1})^{-1}   }{\longrightarrow}} & \bQ ((\tf _0^{1/\infty}))   &  { \overset{\pi  ^{-1}} {\longrightarrow} } &
 \bQ ((\zfk  ^{1/\infty})) &  { \overset{\si  _{\01}} {\longrightarrow} } &  \bQ ((\zfk  ^{1/\infty})) &  { \overset{\pi } {\longrightarrow} } &   \bQ ((\tf _0^{1/\infty})) \, .  \\ \\
 
\end{matrix}
$$

Hence it follows that
$$
d^{-1  }\cdot \si _ {\overrightarrow {   \infty 1  } }\cdot d\cdot (\si _ {Q_2})^{-1  }=k_*(\pi \cdot \si _{\01}\cdot \pi ^{-1}\cdot \si _{Q_1}^{-1}))=
k_*(\pi \cdot \si _{\01}\cdot \pi ^{-1}\cdot \si _{Q_1}^{-1} \cdot \pi \cdot \pi ^{-1})=
$$
$$
k_*(\pi \cdot (\ffk _{\pi} (\si ) )^{-1}\cdot \pi ^{-1})=d^{-1}\cdot (k_* (\ffk _{\pi} (\si )))^{-1}\cdot d\, .
$$

\noindent
Therefore  $q_1^ {-1  }\cdot(d^{-1  }\cdot \si _ {\overrightarrow {   \infty 1  } }\cdot d\cdot (\si _ {Q_2})^{-1  })\cdot q_1= q_2 ^{-1}\cdot k_*\big( \ffk _\pi (\si )\big)^{-1} \cdot q_2$.

\medskip

Now we shall prove the  point ii). Observe that $R_k^*(\uf _k)=\uf _0$ and $R_k^*(\sf _k)=\sf _0$.
We have the following commutative diagram of field isomorphisms

$$
\begin{matrix}

 \bQ ((\uf _\sb ^{1/\infty})) &{\overset{(\si _ {  \overrightarrow{\infty \xn ^1}      })^{-1} }{\longrightarrow}} &     \bQ ((\uf _1^{1/\infty})) &{ \overset{  c_1} {\longrightarrow} }&  \bQ ((\sf _1^{1/\infty})) &  { \overset{    \si    _{     Q_3(1)    }                    } {\longrightarrow} } &    \bQ ((\sf _\sb ^{1/\infty})) &   { \overset{ c_\sb ^{-1}   } {\longrightarrow} }& \bQ ((\uf _ \sb ^{1/\infty}))  \\ \\
 {R_\sb ^*  }\downarrow &&{R_1^* }\downarrow    && {R_1 ^* }\downarrow    && {R_\sb^* }\downarrow  &&  {R_\sb ^* }\downarrow  \\ \\
 \bQ ((\uf _0^{1/\infty}))    &{\overset{ (\si _ {  \overrightarrow{\infty 1}    })^{-1}   }{\longrightarrow}} & \bQ ((\uf _0^{1/\infty}))   &  { \overset{k(\pi )  } {\longrightarrow} } &
 \bQ ((\sf _0  ^{1/\infty})) &  { \overset{\si  _{Q_2}} {\longrightarrow} } &  \bQ ((\sf _0  ^{1/\infty})) &  { \overset{k(\pi )^{-1} } {\longrightarrow} } &   \bQ ((\uf _0^{1/\infty})) \, .  \\ \\
 
\end{matrix}
$$
Hence it follows that $(R_\sb )_* (k(\pi )^{-1} \cdot \si _{Q_2}\cdot k(\pi ) \cdot ( \si _  {\overrightarrow { \infty 1     }}  )^{-1})=
c_{\sib} ^{-1 }\cdot  \si _ {Q_3(1) }\cdot  c_1\cdot ( \si _ {\overrightarrow { 0 {\xn ^1}     } })^{ -1}
$.
Observe that $k(\pi )^{-1} \cdot \si _{Q_2}\cdot k(\pi ) \cdot ( \si _  {\overrightarrow { \infty 1     }}  )^{-1}=k_*(\ffk _\pi (\si ))$. Hence
$c_{\sib} ^{-1 }\cdot  \si _ {Q_3(1) }\cdot  c_1\cdot ( \si _ {\overrightarrow { 0 {\xn ^1}     } })^{ -1}=(R_\sb)_*(k_* (\ffk _\pi (\si )))$. 

The proof of the point iii) is similar and the point iv) is obvious.
\hpb

\bigskip

\noindent
{\bf Lemma 2.6.}
{\it Let $\si \in G_\Qbb$. Then in the group $\pi _1^{\rm et}(V,\01 )$ we have
\begin{enumerate}
\item[i)] $ \pi ^{-1}\cdot \big( r^{-1}\cdot  \si _ {\jn     \overrightarrow {1 \infty    }          }\cdot r \cdot ( \si _ {\jn \overrightarrow {1 0     } })^{-1}\big) \cdot \pi=y_0^{\frac{1-\chi (\si )}{2} }$,

\medskip

\item[ii)] $q_2^{-1}\cdot \big( t_{\sib} ^{-1 }\cdot  \si _ {\overrightarrow {\infty  {\xn ^1}      } }\cdot  t_1\cdot ( \si _ {\overrightarrow { \infty  1     } })^{ -1})\cdot  q_2= z^{\frac{\sib -\chi (\si )}{N}  }$,

\medskip

\item[iii)] $ q_4(\sib)^{-1 }\cdot \big( u_{\sib} ^{-1 }\cdot  \si _ {Q_4(1)    }\cdot  u_1\cdot ( \si _ {Q_3(1)})^{ -1}\big)\cdot
q_4(\sib)= (y_   {\sib  -1}\cdot \ldots \cdot y_2\cdot y_1)^{-1}\cdot  y_ {\sib }^{\frac{1-\chi (\si )}{2} }   \cdot (    y_   {\sib  -1}\cdot \ldots \cdot y_2\cdot y_1)$,

\medskip

\item[iv)] $ q_6(\sib)^{-1 }\cdot \big( (s_{ \sib })^{-1}\cdot \si _{\01 }\cdot s_1\cdot  ( \si _ {\overrightarrow {0  {\xn ^1}  }   })^{ -1}\big) \cdot q_6(\sib )= (y_   {\sib  -1}\cdot \ldots \cdot y_2\cdot y_1)^{-1}\cdot x^{\frac{\sib -\chi (\si  )}{N} }   \cdot (    y_   {\sib  -1}\cdot \ldots \cdot y_2\cdot y_1)$.

 \end{enumerate}
}

\smallskip

\noindent
{\bf  Proof.}  We shall show the point ii). The remaining identities can be proved in the same way.
We shall act by the isomorphism $ t_{\sib} ^{-1 }\cdot  \si _ {\overrightarrow {\infty  {\xn ^1}      } }\cdot  t_1\cdot ( \si _ {\overrightarrow { \infty  1     } })^{ -1}$
on $\uf _0^{1/n}$.
We have
$$
\uf _0^{1/n}\overset{ ( \si _ {\overrightarrow { \infty  1     } })^{ -1}}{\longmapsto}
\uf _0^{1/n}\overset{t_1}{\longmapsto} \xn ^{-1/n}   \uf _1^{1/n}    
 \overset{\si _ {\overrightarrow{\infty { \xn ^1 } }}   } {\longmapsto}  \xi _{Nn}^{-\chs}\uf_\sb ^{1/n}\overset{t_\sb ^{-1}}{\longmapsto}\xi _{Nn}^{-\chs}\xi _{Nn}^{\sb}\uf _0^{1/n}
=\xi _n ^{\frac{\sb - \chs}{N}}\uf _0^{1/n}.
$$
Therefore $ t_{\sib} ^{-1 }\cdot  \si _ {\overrightarrow {\infty  {\xn ^1}      } }\cdot  t_1\cdot ( \si _ {\overrightarrow { \infty  1     } })^{ -1}=z_\infty  ^{\frac{\sb - \chs}{N}}\uf _0^{1/n}$ and
$q_2^{-1}\cdot z_\infty  ^{\frac{\sb - \chs}{N}}\uf _0^{1/n}   \cdot q_2=z  ^{\frac{\sb - \chs}{N}}\uf _0^{1/n}$.
 \hpb
\bigskip

\noindent
{\bf Theorem 2.7.} (Octagonal relation)
{\it  Let $\si \in G_\Qbb$. Then we have the following identity  in the group $\pi _1 ^{\rm et}(V,\01 )$
$$
1= y_   {\sib  -1}\cdot \ldots \cdot y_2\cdot y_1\, \cdot \,  (y_   {\sib  -1}\cdot \ldots \cdot y_2\cdot y_1)^{-1}\cdot x^{\frac{\sib -\chi (\si  )}{N} }   \cdot (    y_   {\sib  -1}\cdot \ldots \cdot y_2\cdot y_1)\, \cdot
$$
$$
   q_6(\sib)^{-1 }\cdot \big( (R_{\sib})_*\big(  \ffk _\pi (\si )\big)^{-1  } \big) \cdot q_6(\sib)             \, \cdot \,  (y_   {\sib  -1}\cdot \ldots \cdot y_2\cdot y_1)^{-1}\cdot  y_ {\sib }^{\frac{1-\chi (\si )}{2} }   \cdot (    y_   {\sib  -1}\cdot \ldots \cdot y_2\cdot y_1) \, \cdot 
$$
$$
 q_3(\sib)^{-1 }\cdot \big( (R_{\sib})_*\big( k_*\big( \ffk _\pi (\si )\big)\big) \big) \cdot q_3(\sib)    \, \cdot \,  z^{\frac{\sib -\chi (\si )}{N}  }
\, \cdot \, q_2^{-1  }\cdot k_*\big( \ffk _\pi (\si )\big)^{-1} \cdot q_2      \,  \cdot \,  y_0^{\frac{1-\chi (\si )}{2} }  \, \cdot \,\ffk _\pi(\si )\, .
$$
}

\smallskip

\noindent
{\bf  Proof.}  The theorem follows   from  Lemmas 2.1,  2.4, 2.5  and 2.6.
\hpb

\medskip

\noindent
Remark. The special cases of Theorem 2.7 are proved in \cite[Proposition 2.1 and Proposition 2.4]{W10}.

\bigskip

\section{Measures }
\smallskip

  In this section and the next ones we change slightly the notation introduced in the    section 2. We denote by $\Nbb$ the set of non-negative integers. 
Let us fix a rational prime $p$.
 For $n\in \Nbb$  let $V_n:=\Pbb ^1_ {\overline \Qbb }\setminus (\{ 0,\infty \}\cup \mu _{{p^n}})$. Then $V_0=\Pbb ^1_ {\overline \Qbb }\setminus \{ 0,1,\infty \} $.
We denote by $\pi _1(V_n,\01 )$ the maximal  pro-$p$ quotient of the \'etale fundamental group of $V_n$.
We denote by $x_n$ (loop around $0$) and $y_{i,n}$  (loop around ${\xi _{p^n}^i}$) for $i \in \Zbb /(p^n)$  and  $z_n$ (loop around $\infty$)
the standard generators of $\pi _1(V_n,\01 )$
constructed in section 2 if $N=p^n$. However because of certain asymmetry we shall take also indices between $0$ and $p^n$. Let  $\pi _n : \01 \leadsto \frac{1}{p^n} \10$ be  the canonical path along the real  interval $[0,1]$.

Let 
$
 \Yc _n
$
be a set of non-commuting variables $X_n$ and $Y_{k,n}$ for $k\in \bZ/p^n\bZ$ 
and let 
$
 \Qbb  _p\{\{\Yc _n\}\}
$
be a $\Qbb  _p$-algebra of non-commutative formal power series on elements of $\Yc_n$. We denote by $ \Qbb  _p\{\{\Yc _n ^0\}\}$ a $\Qbb  _p$-algebra of non-commutative formal power series on elements $Y_{k,n}$ for $k\in \bZ/p^n\bZ$.
We denote by
$
I_n
$
the augmentation ideal of $ \Qbb  _p\{\{\Yc _n\}\}$.
Let
\[
 E_n:\pi _1(V_n,\01 )\to \Qbb  _p\{\{\Yc _ n\}\}
\]
be a continuous, multiplicative map defined by
\[
 E_n(x_n):=\exp X_n\;\;{\rm and}\;\;E_n(y_{k,n}):=\exp Y_{k,n}\;\;{\rm for}\;\; k\in \bZ/p^n\bZ\,.
\] 
It follows from the Baker-Campbell-Hausdorff formula (see \cite[Theorem 5.19]{MKS}) that for any $g\in \pi _1(V_n,\01 )$, $E_n(g)=\exp G$ for some $G\in  \Qbb  _p\{\{\Yc _ n\}\}$, which is a (possibly infinite) sum of homogeneous Lie elements of  $\Qbb  _p\{\{\Yc _ n\}\}$. Further such a Lie series we shall call a Lie element.
Let $\si \in G_\Qbb$. We recall that $\ffk _{\pi _n}(\si ) = \ffk _{\pi _n}(\si )(x_n,y_{0,n},y_{1,n},\ldots ,y_{p^n-1,n}):=\pi _n ^{-1}\cdot \si (\pi _n)$.

Let $\Mc _n$ be the set of all monomials of positive degrees  in non-commutative variables  $X_n$ and $Y_{k,n}$ for $k\in \bZ/p^n\bZ$ and let $\Mc _n^i$ be the subset of $\Mc _n$ of all monomials of degree $i$. Then for
 $\si \in G_\Qbb$,
\[
E_n(\pi_n ^{-1}\cdot \si (\pi _n))=1+\sum_{w\in \Mc _n}\lam _w^{(n)  }(\si )w =1+\sum _{i=1}^\infty \sum_{w\in \Mc _n^i}\lam _w^{(n)  }(\si )w
\]
for some coefficients $\lam _w^{(n)  }(\si ) \in \Qbb _p$.
The covering map 
\[
pr_n^{n+m}:V_{n+m}\to V_n,\;  pr_n^{n+m}(\zfk )=\zfk ^{p^m}
\]
induces the map of $\Qbb _p$-algebras
\[
\Qbb _p\{\{  \Yc _{n+m}    \}\} {\overset { (pr_n^{n+m})_*} { \longrightarrow} }  \Qbb _p\{\{ \Yc _n  \}\}
\]
such that  $(pr_n^{n+m})_*(X_{n+m}) =p^m X_n$,

\noindent
$   ( pr_n^{n+m})_* (Y_{i+kp^n,n+m})=\exp  (-kX_n)\cdot Y_{i,n}\cdot \exp (kX_n)$ for $0\leq i<p^n$ and $0\leq k<p^m$
and
$$
(pr_n^{n+m})_*\big( E_{n+m}(\pi_{n+m} ^{-1}\cdot \si (\pi _{n+m}))\big) =E_n(\pi_n ^{-1}\cdot \si (\pi _n))\, .
$$

\medskip

\noindent
Let $r\geq 1$.
Let ${\bold i}=(i_1,\ldots ,i_r)\in (\Zbb/(p^n))^r$ be a multi-index. We set $\bold Y_{\bold i}:=Y_{ i_1,n}Y_{ i_2,n}\ldots Y_{ i_r,n}$.
Then   the family of functions 
\[
\big\{  (\Zbb/(p^n))^r\ni {\bold i}\mapsto \lam _{\bold Y_{\bold i}   }^{(n)  }(\si )   \big\} _{n\in \Nbb   }
\]
form a measure on  $(\Zbb _p)^r$ with values in $\Qbb _p$, which we denote by $\Kc _r(\si )$ (see \cite[Proposition 2.4.]{W11} and also \cite[Proposition 6]{NW3}).
We correct below the result stated in \cite[Proposition 1.1. iii)]{W10}.  

\smallskip

\noindent 
{\bf Theorem 3.1.} {\it  Let $0\leq i_1,\ldots ,i_r<p^n$ and let $w=$

\noindent
$ X_n^{n_  0 }Y_{i_1,n}X_n^{n_ 1  }Y_{i_2,n}\ldots X_n^{n_{r-1}   }Y_{i_r,n}X_n^{n_ r  }$. Let ${\bold i  }=(i_1,\ldots ,i_r)$.
Let $a(k) =x_k - x_{k+1}-i_k +i_{k+1}$ for $1\leq k \leq r-1$.
Then for  $\si \in G_\Qbb$,

$$\lam _w^{(n)  }(\si )=\Big( \prod _{i=0}^r n_i!\Big) ^{-1}$$
$$ \int _{{\bold i}+p^n(\Zbb _p)^r}   \Big(  {\frac{i_1-x_1}{p^n}    }\Big)^{n_0} 
 \Big(  {\frac{a(1)}{p^n}    }\Big)^{n_1}...       \Big(  {\frac{a(r-1)}{p^n}    }\Big)^{n_{r-1}}  \Big(  {\frac{x_r-i_r}{p^n}    }\Big)^{n_r}
d\Kc _r(\si )(x_1,\ldots ,x_r)\,. $$ }     

\medskip

\noindent
{\bf  Proof.}  Let $0\leq \al _i <p^m$ for $1\leq i \leq r$. Let us set ${\mathbf \alpha} :=(\al _1,\ldots ,\al _r )$. We have

\smallskip

\noindent
$pr _n^{n+m}(Y_{i_1+\al _1p^n,n+m}Y_{i_2+\al _2p^n,n+m}\ldots Y_{i_r+\al _rp^n,n+m})=$

\smallskip

\noindent
$e^{-\al _1X_n}Y_{i_1,n}e^{\al _1X_n}e^{-\al _2X_n}Y_{i_2,n}e^{\al _2X_n}\ldots e^{-\al _rX_n}Y_{i_r,n}e^{\al _rX_n}$.
Observe that 

\noindent
$pr _n^{n+m}(X_{n+m})=p^mX_n$.
Hence it follows that
$$
\lam _w^{(n)}(\si )\equiv
$$
$$
 \Big( \prod _{i=0}^r n_i!\Big) ^{-1}    \sum _{{\bf \al}=(\al _1,\ldots ,\al _r )\in (\Zbb /(p^n))^r}(-\al _1)^{n_0}(\al _1-\al _2)^{n_1}\ldots (\al _{r-1}-\al _r)^{n_{r-1}}\al _r^{n_r} \lam _{{\bold i}+{\bold \al}p^n}^{(n+m)}(\si ) 
$$
modulo $p^m$.
Observe that the right hand side tends to
$
{\frac{1}{n_0!n_1!\ldots n_r!}}$
$$
 \int _{{\bold i}+p^n(\Zbb _p)^r}   \Big(  {\frac{i_1-x_1}{p^n}    }\Big)^{n_0} 
 \Big(  {\frac{a(1)}{p^n}    }\Big)^{n_1}...       \Big(  {\frac{a(r-1 )}{p^n}    }\Big)^{n_{r-1}}  \Big(  {\frac{x_r-i_r}{p^n}    }\Big)^{n_r}
d\Kc _r (\si ) (x_1,\ldots ,x_r)\, .$$   \hpb

\bigskip

Let $c\in \zp$. We denote by  $\langle c \rangle _n$ an integer such that $0\leq \langle c \rangle _n <p^n$ and $c\equiv \langle c \rangle _n$
modulo $p^n$.  Let $\si \in G_\Qbb$. We set  $\langle \si \rangle _n := \langle \chi (\si ) \rangle _n$. If $n$ is fixed we shall also write  $\langle c \rangle $ 
instead of 
$\langle c \rangle _n$  and  $\langle \si \rangle $ instead  of $\langle \si \rangle _n$.

\medskip

Observe that
\[
\lam _{X_n^\mu Y_{i,n }X_n^\nu}^{(n)  }(\si )={\frac{1}{\mu !\nu !}}\int _{i+p^n\Zbb _p}(-1)^\mu \Big( {\frac{x-i}{p^n}}\Big) ^{\mu +\nu }d\Kc _1(\si )(x)\; .
\]
The coefficient $\lam _{ X_n Y_{i,n } }^{(n)  }(\si )$ is a kind of $\ell$-adic Galois dilologarithm  evaluated at $\xi _{p^n}^{-i}$.
Below we shall show an inversion formula for these coefficients. To simplify the notation we set
\[
\lam _i^{(n)  }(\si ):=\lam _{  Y_{i,n } }^{(n)  }(\si )\,.
\]

\smallskip

\noindent
{\bf Theorem 3.2. } {\it Let $\mu >0$ and $0<i<p^n$ . Then
\[
\int _{i+p^n\Zbb _p} \Big( {\frac{x-i}{p^n}}\Big) ^{\mu   }d\Kc _1(\si )(x) +(-1)^{ \mu +1 }  \int _{p^n -i+p^n\Zbb _p}  \Big( {\frac{x-(p^n-i)}{p^n}}\Big) ^{\mu  }d\Kc _1(\si )(x)=
\]
\[
\sum _{j=0}^{\mu -1}{\mu \choose j  } \int _{i+p^n\Zbb _p} \Big( {\frac{x-i}{p^n}}\Big) ^{j }d\Kc _1(\si )(x)+
\]
\[
{\frac{(-1)^\mu }{p^{n\mu} } }\sum _{j=0}^{\mu }{\mu \choose j  }(i-p^n)^{\mu -j}{\frac { p^{nj}}{j+1 } } \Big(  B_{j+1}\big( {\frac{p^n-i}{p^n}}\big)-
\chi (\si )^{j+1 }  B_{j+1}\big( {\frac{\langle  \chi (\si )^{-1}(   p^n-i)\rangle }{p^n}}\big)  \Big)\;.
\]
}

\medskip

\noindent {\bf Proof.}  Let $0<i<p^n$. Observe that
\[
\int _{i+p^n\Zbb _p} \Big( {\frac {x-i}{p^n}}\Big) ^{\mu   }d\Kc _1(\si )(x)\equiv \sum _{k=0 }^{p^m-1}k^\mu \lam _{i+p^nk}^{(n+m)}(\si )\;{\rm mod}\; (p^m)\, .
\]
We have
\[
 \sum _{k=0 }^{p^m-1}k^\mu \lam _{p^n-i+p^nk}^{(n+m)}(\si )= \sum _{k=0 }^{p^m-1}k^\mu \lam _{p^{n+m}-(p^n-i+p^nk)}^{(n+m)}(\si )+ \sum _{k=0 }^{p^m-1}k^\mu E_{1,\chi (\si )}^{(n+m)}(p^n-i+p^nk)
\]
(see \cite[Lemma 4.1.]{W11}).
After  calculations we get that
\[
 \sum _{k=0 }^{p^m-1}k^\mu \lam _{p^{n+m}-(p^n-i+p^nk)}^{(n+m)}(\si )\equiv (-1)^\mu \sum _{j=0}^\mu {\mu \choose j}\int _{i+p^n\Zbb _p} \Big( {\frac{x-i}{p^n}}\Big) ^{j   }d\Kc _1(\si )(x) \;{\rm mod}\; (p^m)\,
\]
and
\[
\sum _{k=0 }^{p^m-1}k^\mu E_{1,\chi (\si )}^{(n+m)}(p^n-i+p^nk)\equiv 
\]
\[
{\frac{1 }{p^{n\mu} } }\sum _{j=0}^{\mu }{\mu \choose j  }(i-p^n)^{\mu -j}{\frac { p^{nj}}{j+1 } } \Big(  B_{j+1}\big( {\frac{p^n-i}{p^n}}\big)-
\chi (\si )^{j+1 }  B_{j+1}\big( {\frac{\langle  \chi (\si )^{-1}(   p^n-i)\rangle }{p^n}}\big)  \Big)\; 
\]
 {\rm mod}  $ (p^m)$ .
\hpb

\medskip

\noindent
{\bf Corollary 3.3. } {\it Let $0<i<p^n$. We have
\[
\int _{i+p^n\Zbb _p} \Big( {\frac{x-i}{p^n}}\Big) d\Kc _1(\si )(x) +  \int _{p^n-i+p^n\Zbb _p} \Big( {\frac{x-(p^n-i)}{p^n}}\Big) d\Kc _1(\si )(x) =
\]
\[
-\lam _{  Y_{i,n } }^{(n)  }(\si )+{\frac {(i-p^n)}{p^n}}\Big( B_1\big(   {\frac {(p^n-i)}{p^n}}\big) -\chi (\si )B_1\big(   {\frac {\langle \chi (\si )^{-1}(p^n-i)\rangle }{p^n}}\big)\Big) +
\]
\[
   {\frac {1}{2}} \Big( B_2\big(   {\frac {(p^n-i)}{p^n}} \big) -\chi (\si )^2B_2\big(   {\frac {\langle \chi (\si )^{-1}(p^n-i)\rangle }{p^n}}\big) \Big) \;=
\]
\[
-\lam _{  Y_{p^n-i,n } }^{(n)  }(\si )-{\frac {i}{p^n}}\Big( B_1\big(   {\frac {i}{p^n}}\big) -\chi (\si )B_1\big(   {\frac {\langle \chi (\si )^{-1}i\rangle }{p^n}}\big)\Big) +
\]
\[
   {\frac {1}{2}} \Big( B_2\big(   {\frac {i}{p^n}} \big) -\chi (\si )^2B_2\big(   {\frac {\langle \chi (\si )^{-1}i\rangle }{p^n}}\big) \Big) \;.
\]

} \hpb

\medskip

\noindent
Remark. The coefficient  $\lam _{X_n^\mu Y_{i,n }X_n^\nu}^{(n)  }(\si )$ is a kind of $\ell$-adic Galois polylogarithm evaluated at $\xi _{p_n}^ {-i}$. It would be interesting to compare the formulas from Theorem 3.2 and Corollary 3.3 with the inversion formula for $\ell$-adic Galois polylogarithms in \cite[formula 6.31]{NW2}.

\bigskip

\section{Octagonal relation modulo the third power of  the augmentation ideal }

\medskip

In this section $n$ is fixed. We recall from section 3 that the map $E_n:\pi _1(V_n,\01 )\to \Qbb _p \{\{\Yc _n\}\}$ is continuous and multiplicative.

\medskip

\noindent
{\bf  Lemma  4.1. } {\it  Let $\si \in \Gal (\overline \Qbb / \Qbb)$. Then we have the following identity between formal power series in $ \Qbb _p \{\{\Yc _n\}\}$

$$
1=E_n\big( y_   {\sibn  -1,n}\cdot \ldots \cdot y_{1,n} \big)\, \cdot \,E_n\big(  (y_   {\sibn  -1,n}\cdot \ldots \cdot y_{1,n})^{-1}\cdot x_n^{\frac{\sibn -\chi (\si  )}{p^n} }   \cdot (    y_   {\sibn  -1,n}\cdot \ldots \cdot y_{1,n})\big)\, \cdot
$$
$$
 E_n\big(  q_6(\sibn)^{-1 }\cdot \big( (R_{\sibn})_*\big(  \ffk _{\pi _n} (\si )\big)^{-1  } \big) \cdot q_6(\sibn)     \big)        \, \cdot
$$
$$
 E_n\big(  (y_   {\sibn  -1,n}\cdot \ldots \cdot  y_{1,n})^{-1}\cdot ( y_ {\sibn ,n })^{\frac{1-\chi (\si )}{2} }   \cdot (    y_   {\sibn  -1,n}\cdot \ldots \cdot  y_{1,n})\big) \, \cdot 
$$
$$
E_n\big( q_3(\sibn)^{-1 }\cdot \big( (R_{\sibn})_*\big( k_*\big( \ffk _{\pi _n} (\si )\big)\big) \big) \cdot q_3(\sibn)  \big)  \, \cdot \,E_n\big(  z_n^{\frac{\sibn -\chi (\si )}{p^n}  }\big)
\, \cdot \,
$$
$$
 E_n\big( q_2^{-1  }\cdot k_* \big( \ffk _{\pi _n} (\si )\big)^{-1} \cdot q_2  \big)    \,  \cdot \,
E_n\big(  y_{0,n}^{\frac{1-\chi (\si )}{2} } \big) \, \cdot \,E_n\big(\ffk _{\pi _n}(\si )\big)\, .
$$
}

\medskip

\noindent {\bf Proof.} The lemma follows by applying the map $E_n$ to the identity of Theorem 2.7.   \hpb

\medskip

In order to make more clear the notation let us set
\[
\al _i^{(n)}(\si ):=\lam _{  Y_{i,n } }^{(n)  }(\si )\,,\;\; \be _{a,b}^{(n)}(\si ):=\lam _{  Y_{a,n }Y_{b,n} }^{(n)  }(\si )\,,\;\;\ga _i^{(n)}(\si ):=\lam _{X_n  Y_{i,n } }^{(n)  }(\si )\,.
\]
We shall also omit the subscript $n$ in $Y_{i,n}$ as $n$ is fixed. The summation $\sa$ means $\sum _{(a,b)\in ( \Zbb /(p^n))^2}$.

\medskip

Let us set
\[
A_n:= 1+ \sum _{i=0}^{p^n-1} \al _i ^{(n)} (\si )Y_ i+\sa  \be _{a,b}^{(n)}   (\si )Y_aY_b +\sip \ga _i^{(n)}(\si )\big( X_nY_i-Y_iX_n\big) \, ,
\]

\medskip

\[
B_n:= 1+{\frac{1-\chi (\si )}{2}}Y_0+ {\frac{(1-\chi (\si ) )^2}{8}}Y_0 ^2\, ,
\]

\medskip

\[
C_n:= 1-   \sum _{i=0}^{p^n-1} \al _{p^n-i} ^{(n)} (\si )Y_ i-  \sum _{i=0}^{p^n-1} \al _{p^n-i} ^{(n)} (\si ) X_n Y_i+  \sum _{i=0}^{p^n-1} \al _{p^n-i} ^{(n)} (\si )Y_  iX_n +
\]
\[
  -    \sum _{b=1}^{p^n-1} \al _{p^n-b} ^{(n)} (\si )  \big( \sum _{b<a\leq p^n}Y_a\big) Y_b+ \sum _{a=1}^{p^n-1} \al _{p^n-a} ^{(n)} (\si )Y_a \big( \sum _{a<b\leq p^n}Y_b\big) -\sa \be _{-a,-b }^{(n)} (\si )Y_aY_b +
\]
\[
\sbp \ga _{ -b }^{(n)} (\si ) \big( X_n+ \sap Y_a\big) Y_b-  \sap \ga _{ -a }^{(n)} (\si )Y_a \big( X_n+ \sbp Y_b\big) + \big( \sip \al _{p^n-i} ^{(n)} (\si )Y_ i \big) ^2 \, ,
\]

\medskip

\[
D_n:=   1+ {\frac{\chs - \ls}{ p^n}}X_n + {\frac{\chs - \ls}{ p^n}}(\sip Y_i)+{\frac{\chs - \ls}{2 p^n}}(Y_0X_n-X_nY_0)+
\]
\[
 {\frac{\chs - \ls}{2  p^n}}\big( \sip (X_nY_i  - Y_iX_n)\big) + {\frac{\chs - \ls}{2  p^n}}  \big( \sum _{1\leq b<a\leq p^n}Y_aY_b \big) +
\]
\[
  ( -1)  {\frac{\chs - \ls}{2  p^n}} \big( \sum _{1\leq a<b\leq p^n}Y_aY_b\big)  +  {\frac{1}{2}}\big(  {\frac{\chs - \ls}{  p^n}}\big) ^2 \big( \sa Y_aY_b \big)+
\]
\[
   {\frac{1}{2}}\big(  {\frac{\chs - \ls}{  p^n}}\big) ^2 \big( X_n(\sip Y_i )+(\sip Y_i)X_n\big) \, ,
\]

\medskip

\[
F_n : = 1+\sum _{a=1}^{\ls}\al _{\ls -a}^{(n)}(\si )\Big(  Y_a+Y_a \big( \sum _{j=1}^{a-1}Y_j\big) -\big( \sum _{j=1}^{a-1}Y_j\big) Y_a \Big)+ \al _{\ls }^{(n)}(\si )Y_0 +
\]
\[
 \sum _{ a=\ls +   1}^{p^n-1}\al _{\ls -a}^{(n)}(\si )\Big(  Y_a- Y_a \big(X_n +\sum _{j=a+1}^{p^n}Y_j\big)+  \big(X_n +\sum _{j=a+1}^{p^n}Y_j\big) Y_a\Big) + 
\]
\[
 \sa \be _{\ls -a,\ls -b}^{(n)}(\si )Y_aY_b - \sum _{j=0}^{p^n-1} \ga _{\ls -j}^{(n)}(\si )\Big( \big(X_n +\sum _{k=0}^{p^n-1}Y_k\big)Y_j - Y_j   \big(X_n +\sum _{k=0}^{p^n-1}Y_k\big) \Big)\,  ,
\]

\medskip

$$G_n:=1+   {\frac{1-\chi (\si )}{2}}Y_{\ls }  +             $$
\[
        {\frac{1}{2}}  \big(   {\frac{1-\chi (\si )}{2}} \big) ^2(Y_{\ls })^2+ {\frac{1-\chi (\si )}{2}}\Big( Y_{\ls }( \sum _{i=1}^{\ls - 1}Y_i) - ( \sum _{i=1}^{\ls - 1}Y_i) Y_{\ls }\Big)\, ,
\]

\medskip

$$H_n:=1-       \sip \al _{i- \ls }^{(n)}(\si )\Big(  Y_i +Y_i ( \sum _{j=1}^{\ls -1}Y_j )-  ( \sum _{j=1}^{\ls -1}Y_j )Y_i \Big)             +            $$
\[
   (-1)  \sum _{i=0}^{\ls -1} \al _{i- \ls }^{(n)}(\si ) \big( Y_iX_n -X_nY_i \big)  - \sa \be _{a-\ls ,b- \ls }^{(n)}(\si )Y_aY_b +
\]
\[
  (-1) \sip  \ga _{i-\ls }^{(n)}(\si )\big( X_nY_i - Y_iX_n \big) +  \big( \sip  \al _{i- \ls }^{(n)}(\si ) Y_i \big) ^2 \, ,
\]

\medskip

\[
J_n :=  1+
\]
\[
 {\frac{ \ls - \chs}{ p^n}}X_n + \sum _{i=1}^{\ls - 1}  {\frac{\chs - \ls}{ p^n}}  \big( X_n Y_i - Y_iX_n\big) + {\frac{1}{2}}  \Big(   {\frac{\ls - \chs }{p^n}} \Big) ^2X_n^2\, ,
\]

\medskip

\[
K_n:= 1+ \sum _{i=1}^{\ls - 1}Y_i +  {\frac{1}{2}} \sum _{1\leq j < i< \ls }\big(  Y_iY_j - Y_jY_i \big)+ {\frac{1}{2}} \Big(  \sum _{i=1}^{\ls - 1}Y_i  \Big) ^2\,  .
\]

\medskip

\noindent
Let $\si \in G_\Qbb$. We recall that $\ffk _{\pi _n}(\si ) = \ffk _{\pi _n}(\si )(x_n,y_{0,n},y_{1,n},\ldots ,y_{p^n-1,n})=\pi _n ^{-1}\cdot \si (\pi _n)$.

\medskip

\noindent
{\bf Lemma 4.2}  {\it  Let $\si \in G_ \Qbb$. Let us set $a= y_   {\sibn  -1,n}\cdot \ldots \cdot y_{2,n}\cdot y_{1,n}$ if $\sibn \neq 1$ and $a=1$ if $\sibn =1$.                            We have
\begin{enumerate}
\item[i)] $E_n\big( \ffk _{\pi_n} (\si )(x_n,y_{0,n},\ldots ,y_{p^n-1,n})\big) \equiv A_n\;\;{\rm mod}\;\; I_n^3$

\smallskip

\item[ii)] $E_n\big( ( y_{0,n})^{\frac{1-\chi (\si )}{2}}\big)\equiv B_n \;\;{\rm mod}\;\; I_n^3$

\smallskip

\item[iii)] $E_n\big(  q_2^{-1  }\cdot k_*\big( \ffk _\pi (\si ) (x_n,y_{0,n},\ldots ,y_{p^n-1,n})  \big)^{-1} \cdot q_2  \big)    \equiv C_n\;\;{\rm mod}\;\; 
 I_n^3$

\smallskip

\item[iv)] $E_n\big( ( z_n) ^{\frac{\sibn -\chi (\si )}{p^n}  }\big)  \equiv D_n\;\;{\rm mod}\;\; I_n^3$

\smallskip

\item[v)] $E_n\big(  q_3(\sibn)^{-1 }\cdot \big( (R_{\sibn})_*\big( k_*\big( \ffk _{\pi _n} (\si ) (x_n,\ldots )         \big)\big) \big) \cdot q_3(\sibn) \big) \equiv F_n\;\;{\rm mod}\;\; I_n^3$

\smallskip

\item[vi)] $E_n\big( a^{-1}\cdot  (y_ {\sibn ,n }) ^{\frac{1-\chi (\si )}{2} }   \cdot a \big)  \equiv G_n\;\;{\rm mod}\;\; I_n^3$

\smallskip

\item[vii)] $E_n\big(   q_6(\sibn)^{-1 }\cdot \big( (R_{\sibn})_*\big(  \ffk _{\pi _n} (\si )(x_n,\ldots )       \big)^{-1  } \big) \cdot q_6(\sibn)     \big)  \equiv H_n\;\;{\rm mod}\;\; I_n^3$

\smallskip

\item[viii)] $E_n\big(  a^{-1}\cdot ( x_n)^{\frac{\sibn -\chi (\si  )}{p^n} }   \cdot a\big) \equiv J_n\;\;{\rm mod}\;\; I_n^3$

\smallskip

\item[ix)] $E_n\big( y_   {\sibn  -1,n}\cdot \ldots \cdot y_{2,n}\cdot y_{1,n} \big)  \equiv K_n\;\;{\rm mod}\;\; I_n^3$.
\end{enumerate}
}

\medskip

\noindent
{\bf Proof.} We have $E_n(\ffk _{\pi _n}(\si )) =$
$$
 1+\sum _{w\in \Mc}\lam_w ^ {(n)}(\si )w\equiv 
1+\sip \lam_{Y_{i}}^ {(n)}(\si ) Y_{i}+\sap \sbp \lam_{Y_{a} Y_{b}      }^ {(n)}(\si ) Y_{a} Y_{b}+
$$
$$
\sip \lam_{X_nY_{i}}^ {(n)}(\si )X_n Y_{i} +
\sip \lam_{Y_{i} X_n}^ {(n)}(\si ) Y_{i}X_n  \;\;{\rm modulo}\;\;I_n^3
$$
by the very definition of the coefficients $\lam _w^{(n)}(\si )$ for $w\in \Mc _n$. Observe that
 the coefficient at $X_n$ vanishes. Hence it follows that  
$\lam_{X_nY_{i}}^ {(n)}(\si )+\lam_{Y_{i} X_n}^ {(n)}(\si ) =0$. Therefore $E_n(\ffk _{\pi _n}(\si ))\equiv A_n$ modulo $I_n^3$.
The points ii) -- ix) follow from the point i) and Lemmas 2.1, 2.2, 2.3, 2.5 and 2.6. 
\hpb

\medskip

\noindent
{\bf Proposition 4.3.} {\rm (The pro-$p$ octagonal relation modulo the third powers)}
{\it     Let $\si \in G_\Qbb$. We have the following congruence in $ \Qbb _p \{\{\Yc _n\}\}$
\[
K_n\cdot  J_n \cdot H_n \cdot G_n \cdot F_n \cdot D_n \cdot C_n \cdot B_n \cdot A_n \equiv  1\; \;{\rm modulo}\;\;I_n^3\, .
\]
}

\noindent
{\bf Proof.}  The proposition  follows from Theorem 2.7 and Lemmas 4.1 and 4.2.
\hpb

\medskip

\noindent
{\bf Corollary 4.4.} {\rm Let $\si \in G_\Qbb$. Then 
$$
\sum _{i=0}^{p^n-1} \al _i ^{(n)} (\si )Y_ i+{\frac{1-\chi (\si )}{2}}Y_0 - \sum _{i=0}^{p^n-1} \al _{p^n-i} ^{(n)} (\si )Y_ i+ {\frac{\chs - \ls}{ p^n}}(\sip Y_i)+
$$
$$
\sum _{i=1}^{\ls}\al _{\ls -i}^{(n)}(\si ) Y_i +{\frac{1-\chi (\si )}{2}}Y_{\ls }  -   \sip \al _{i- \ls }^{(n)}(\si ) Y_i + \sum _{i=1}^{\ls - 1}Y_i =0\, .
$$
}

\noindent
{\bf Proof.} The corollary follows immediately
 from the proposition. We write  terms of degree 1 in variables $\Yc _n^0$ starting from the right, i.e. from $A_n$.   \hpb

\medskip

\noindent
Remark. One can easily see that one gets a relation between measures (see Proposition 8.1). It is much more difficult to write terms of degree 2 in variables 
$\Yc _n^0$ and next to arrange them in order to get a relation between measures.

\medskip

The next lemma will be used in section 5, Example 5.2 to show that the coefficients $\ga _i^{(n)} (\si )$ can be adjusted to give a measure on $(\zp )^2$.
 
\medskip

\noindent
{\bf Lemma 4.5.} 
{\it   Let $n$ and $m$ be non negative integers. Let $0\leq i <p^n$ and let $\si \in G_\Qbb$.    Then    we have
\[
\ga _i^{(n)}(\si )=\sum_{k=0}^{p^m-1}p^m \ga _{i+kp^n}^{(n+m)}(\si )+\sum_{k=0}^{p^m-1}(-k)\al _{i+kp^n}^{(n+m)}(\si )
\]
for $0\leq i < p^n$.
}

\medskip

\noindent
{\bf Proof.} The lemma follows from the equality 
$
pr_n^{n+m}\big( E_{n+m}(\pi_{n+m} ^{-1}\cdot \si (\pi _{n+m}))\big) =E_n(\pi_n ^{-1}\cdot \si (\pi _n))\, .
$
\hpb
 
\medskip

\section{Examples of measures }
\smallskip

In this section we present some examples of measures on $\Zbb _p$ and on $(\Zbb _p)^2$. These measures appear when we try to adjust the relations between coefficients in degrees 1 and 2  to get relations between measures. We shall need them in order to write the  formulas of Proposition 8.1 and of Theorems 8.2 and 8.4  in section 8.

\smallskip
\medskip

We denote by $$P: \zp [[\zp ]]\to \zp [[T]]$$  the Iwasawa isomorphism given by the condition $P([1])=1+T$.
For $\mu \in  \zp [[\zp ]]$ we set
$F(\mu )(X):=P(\mu )(\exp X - 1)$. We use the similar definitions and  notations for measures on $(\zp )^m$.

\smallskip

\noindent
Let $c\in \Zbb _p$, $c=\sum _{i=0} ^\infty \al _ip^i$ with $0\leq \al _i <p$.  Observe that  $\langle c\rangle _n:=\sum _{i=0} ^ {n-1} \al _ip^i$.

\medskip

\noindent
{\bf Example  5.1. } {\it Let $c=\sum _{i=0} ^\infty \al _ip^i \in \Zbb _p^\times$ with  $0\leq \al _i <p$ for all $i$. The family of functions
\medskip
\[
\mu ^{(n)}:\Zbb /(p^n)\ni i\mapsto     \left\{
\begin{array}{ll}
{\frac{c-\langle c\rangle _n}{p^n}} +1 & \text{ if } 1\leq i< \langle c \rangle _n, \\
{\frac{c-\langle c\rangle _n}{p^n}} & \text{ if }\langle c \rangle _n \leq i \leq p^n\;.
\end{array}
\right.
\]
\medskip
is a measure on $\Zbb _p$.
We denote this measure by $\Mc (c)$. We have $P(\Mc (c))(T)={\frac{(1+T)^c-(1+T)}{T}}$.
 Observe that 
\[
\Mc (c) =\Big( \sum _{i=0}^{p^n-1}\frac{c- \langle c\rangle _n}{p^n} [i]+\sum _{i=1}^{ \langle c\rangle _n -1}[i] \Big) _{n\in \Nbb} \in \varprojlim _n \zp [\Zbb /(p^n)]\, .
\]
}

\medskip

\noindent {\bf Proof.} 
Let $i$ be such that $1\leq i <  \langle c\rangle _n $. Then $i+kp^n< \langle c\rangle _{n+1} =\langle c\rangle _n +\al _np^n$ if and only if $k\leq \al _n$. Hence we have 
\medskip
\[
\mu ^{(n+1)}:\Zbb /(p^{n+1})\ni i +kp^n\mapsto     \left\{
\begin{array}{ll}
{\frac{c-\langle c\rangle _{n+1}}{p^{n+1}}} +1 & \text{ if } 0\leq k \leq \al _n, \\
{\frac{c-\langle c\rangle _{n+1}}{p^{n+1}}} & \text{ if }\al _n<  k < p\;.
\end{array}
\right.
\]
\medskip
Therefore 
$$\sum _{k=0}^ {p-1}\mu ^{(n+1)}(i+kp^n)=\big( \sum _{k=0}^ {p-1}{\frac{c-\langle c\rangle _{n+1}}{p^{n+1}}}\big) +\al _n+1=
$$
$$
 {\frac{c-\langle c\rangle _n-\al _np^n}{p^n}}+\al _n +1={\frac{c-\langle c\rangle _n}{p^n}} +1=\mu ^{(n)}(i)\, .$$

If $\langle c \rangle _n \leq i \leq p^n$ the proof of the distribution relation is similar.  Let $k>0$. Then the Riemann sum 
$$
  \sip  i^k\cdot {\frac{c-\langle c\rangle _n}{p^n}} +\sum _{i=0}^{ \langle c\rangle _n -1}i^k =
S_k(p^n)\cdot {\frac{c-\langle c\rangle _n}{p^n}}  +S_k( \langle c\rangle _n )\to {\frac{1}{k+1}}\big( B_{k+1}(c) - B_{k+1}\big)
$$
 if $n\to \infty$. Therefore $ \int _{\Zbb _p}x^kd\mu (x)= {\frac{1}{k+1}}\big( B_{k+1}(c) - B_{k+1}\big)$ 
for $k>0$.
For $k=0$, $\int _{\Zbb _p}d\mu (x)=c - 1$.
The power series 
$$F(\mu )(X)=\sum _{k=0}^\infty \big( \int _{\Zbb _p}x^k d\mu (x)\big) {\frac{X^k}{k!}}=$$
$$
c-1+{\frac{1}{X}}\sum _{k=1}^ \infty \big( B_{k+1}(c) - B_{k+1}\big){\frac{X^{k+1}}{(k+1)!}}=-1+{\frac{\exp (cX)-1}{\exp X-1}}= {\frac{\exp (cX)-\exp X}{\exp X-1}}\, . $$ 
Hence it follows that
$P(\Mc (c))(T)={\frac{(1+T)^c-(1+T)}{T}}$.  \hpb

\medskip

\noindent
Remark.  Let $m$ and $p$ be coprime and let $0<a<m$. The measure $\Mc ({\frac{a}{m}})+\de _0$ appears in \cite{W13}, where it is denoted by $\mu ({\frac{a}{m}})$.

\bigskip

\noindent
{\bf Example 5.2. } {\it   Let $\si \in G_\Qbb$.  The family of functions

\medskip
\[
\mu ^{(n)}:\big( \Zbb /(p^n) \big) ^2  \ni (a,b) \mapsto   \ga _{-b}^{(n)}(\si )-\ga _ {-a}^{(n)}(\si )+ \left\{
\begin{array}{ll}
\al _{-a}^{(n)}(\si ) & \text{ if }\;  0<a<b\leq p^n, \\
-\al _{-b}^{(n)}(\si ) & \text{ if }\;  0<b<a\leq p^n,\\
0  & \text{ if }\;  a=b
\end{array}
\right.
\]
\medskip
is a measure on $(\Zbb _p) ^2$, which we denote by $D_2(\si )$. Therefore $D_2(\si )=$
$$
\Big(\sum _{0\leq a,b<p^n}(\ga _{-b}^{(n)}(\si )-\ga _{-a}^{(n)}(\si ))[a,b] +\sum _{0<a<b\leq p^n} \al _{-a}^{(n)}(\si )[a,b]-\sum _{0<b<a\leq p^n}\al _{-b}^{(n)}(\si ) [a,b]\Big) _{n\in \Nbb}\in  
$$
$ \varprojlim _n \zp [\Zbb /(p^n)\times \Zbb /(p^n)   ]$.
}

\bigskip

\noindent
{\bf Proof.}
Let $  0<a<b\leq p^n$. Then $a+\al p^n<b+\be p^n$ if and only if $0\leq \al \leq \be <p$. Then we have 

\smallskip

\noindent
$$\sum _{\al =0}^{p-1} \sum _{\be =0}^{p-1}\mu ^{(n+1)}(a+\al p^n,b+\be p^n)  =     \sum _{0\leq  \al \leq \be <p}  \mu ^{(n+1)}  (a+\al p^n,b+\be p^n)+ $$

\noindent
$$\sum _{0\leq  \be < \al <p}  \mu ^{(n+1)}  (a+\al p^n,b+\be p^n)=
\sum _{0\leq  \al \leq \be <p}\al _{- (a+\al p^n)}^{(n+1)}(\si ) $$

\noindent
$$- \sum _{0\leq  \be < \al <p}\al _{- (b+\be  p^n)}^{(n+1)}(\si ) +
\sum _{\be =0}^{p-1}p\ga  _{- (b+\be  p^n)}^{(n+1)}(\si ) - \sum _{\al =0}^{p-1}p\ga  _{- (a+\al   p^n)}^{(n+1)}(\si )=$$

\noindent
$$
\sum _{\al =0}^{p-1}(p- \al )\al _{- (a+\al p^n)}^{(n+1)}(\si ) - \sum _{\be =0}^{p-1}(p-1 - \be )\al _{- (b+\be p^n)}^{(n+1)}(\si ) +
$$

$$
\sum _{\be =0}^{p-1}p\ga  _{- (b+\be  p^n)}^{(n+1)}(\si ) - \sum _{\al =0}^{p-1}p\ga  _{- (a+\al   p^n)}^{(n+1)}(\si )\,.$$

The indices $-(a+\al p^n)$ and $-(b+\be p^n)$ we write as $p^n -a+(p-1-\al )p^n$ and $p^n-b+(p-1-\be )p^n$.
Observe that $0\leq \al  <p$ implies $0\leq p-1-\al <p$. Hence we get
$$
\sum _{\al =0}^{p-1}(\al +1) \al _{p^n - a+\al p^n}^{(n+1)}(\si ) -\sum _{\be =0}^{p-1} \be \al _{p^n- b+\be p^n}^{(n+1)}(\si ) +\sum _{\be =0}^{p-1}p\ga  _{p^n- b+\be  p^n}^{(n+1)}(\si )
$$
$$
 - \sum _{\al =0}^{p-1}p\ga  _{p^n- a+\al   p^n}^{(n+1)}(\si )
=\ga _{-b}^{(n)}(\si ) - \ga _{-a}^{(n)}(\si )+ \al _{-a}^{(n)}(\si )=\mu ^{(n)}(a,b)\;$$
by Lemma 4.3 and the equality $\sum _{\al =0}^{p-1}(\al +1) \al _{p^n - a+\al p^n}^{(n+1)}(\si ) =\al _{-a}^{(n)}(\si )$.
The cases  $0<b<a\leq p^n$ and $0<a=b\leq p^n$ we check in the same way.  \hpb

\medskip

\noindent
Remark. 
 The coefficients $\ga _i^{(n)}$'s are a kind of $p$-adic Galois dilogarithms at $\xi ^i_{p^n}$. They satisfy a kind of distribution formula (see Lemma 4.5.). (More about the distribution formulas for $l$-adic Galois polylogarithms one can find in \cite{NW3}). However they do not give measures on $\zp$. In the above example we have  seen that 
we can construct  a measure on $(\zp )^2$ from the  $\ga _i^{(n)}(\si )$'s.

\medskip

\noindent
{\bf Example 5.3. } {\it  Let $c\in \zp ^\times$.  The family of functions $\mu ^{(n)}:\Zbb /(p^n) \times \Zbb /(p^n) \to \Zbb _p$ given by

\medskip
\[
\mu ^{(n)} (a,b):= \left\{
\begin{array}{ll}
-{\frac{c -\langle c \rangle _n}{p^n}} & \text{ if }\;  1\leq a<b< \langle c \rangle _n  \; \text{or} \; \langle c \rangle _n \leq a<b\leq p^n\, , \\
{\frac{c - \langle c \rangle _n }{p^n}} & \text{ if }\;  1\leq b<a< \langle c \rangle _n  \; \text{or} \;  \langle c \rangle _n  \leq b<a \leq p^n\, ,  \\
0  & \text{otherwise }\, 
\end{array}
\right.
\]
\medskip

\[
+ \left\{
\begin{array}{ll}
-1 & \text{ if }\;  1\leq a<b< \langle c \rangle _n   \, , \\
1 & \text{ if }\;  1\leq b<a< \langle c \rangle _n \,  ,  \\
0  & \text{otherwise }\, 
\end{array}
\right.
\]
\medskip
is a measure on $\Zbb _p\times \Zbb _p$. We shall denote the above measure on $ \Zbb _p \times \Zbb _p $ by $\Nc _2(c)$.
}

\medskip

\noindent
{\bf Proof.}    Let $c =\sum _{i=0}^{\infty }\si _i  p^i$, where $0\leq \si _i < p$. Then ${\langle c \rangle _{n+1}}=\cs +\si _n p^n$.
 Let $1\leq a<b<\cs$ and let $0\leq \al ,\be \leq p-1$.

\begin{enumerate}
\item[i)]  Then $a+\al p^n<b+\be p^n< \cs +\si _n p^n  $ iff $0\leq \al \leq \be\leq  \si _n$. The number of such pairs $(\al ,\be )$ is $\frac{(\si _n +1)(\si _n +2)}{2}$.
\item[ii)]  The number of pairs $(\al ,\be )$ such that $\cs +\si _n p^n \leq a+\al p^n<b+\be p^n\leq p^{n+1}$ is $\frac{(p-\si _n-1)(p-\si _n)}{2}$.
\item[iii)]  The number of pairs $(\al ,\be )$ such that $b+\be p^n<a+\al p^n<\cs +\si _np^n$  is  $\frac{\si _n (\si _n +1)}{2}$.
\item[iv)]  The number of pairs $(\al ,\be )$ such that $\cs +\si _n p^n \leq b+\be  p^n<a+\al p^n\leq p^{n+1}$ is $\frac{(p-\si _n-1)(p-\si _n-2)}{2}$.
 \end{enumerate}
Now one easily checks that $\sum _{\al ,\be}\mu ^{(n+1)}(a+\al p^n,b+\be p^n)=\mu ^{(n)}(a,b)$.
In the same way one shows the distribution relation in all other cases.  \hpb

 \medskip

\noindent
{\bf Example 5.4. } {\it Let ${\bf a}=(a_1,\ldots ,a_r)\in ( \zp )^r$. We denote by $
\de _{\bf a}$ the Dirac measure on $(\zp )^r$. Observe that $P(\de _{\bf a})(T_1,\ldots ,T_r)=\prod _{i=1}^r(1+T_i)^{a_i}$.} \hpb

\bigskip

\section{Operations on  measures}

In this section we present some elementary operations on measures. These operations will be useful to write down the formula expressing symmetries of  measures $\Kc _1(\si )$ and  $\Kc _2(\si )$ in section 8.

Let $\mu \in \zp [[(\zp )^i]]$  and let ${\bf c}=(c_1,\ldots ,c_i)\in (\zp )^i$. Then the family of functions 
\[
\Big( (\zpn )^i \ni (a_1,\ldots ,a_i) \mapsto \mu ^{(n)} (a_1-c_1,\ldots ,a_i-c_i)\in \zp \Big) _{n\in \Nbb}
\]
is a measure on $(\zp )^i$, which we denote by $T_{\bf c}(\mu )$ and call a translation of $\mu $ by $\bf c$.

\medskip

\noindent 
{\bf Lemma 6.1.} {\it We have \, 
\[
\int _{(\zp )^i}x_1^{n_1}x_2^{n_2}\ldots x_i^{n_i}d(T_{\bf c}(\mu ))(x_1,x_2\ldots x_i)=
\]
\[
\int _{(\zp )^i}(x_1+c_1)^{n_1}(x_2+c_2)^{n_2}\ldots ( x_i+c_i)^{n_i}d\mu (x_1,x_2\ldots x_i)
\]
and 
\[
P( T_{\bf c}(\mu ))(T_1,T_2\ldots T_i)=\, P(\mu )(T_1,T_2\ldots T_i)\prod _{k=1}^i (1+T_k)^{c_k}\, .
\]
}

\medskip

\noindent
{\bf Proof.} One verifies that the formulas hold for $\mu \in \zp [( \zp )^i ]$. Hence by continuity they hold for $\mu \in \zp [[(\zp )^i ]]$. \hpb

\medskip

 Let $t_{\bf c} :(\zp )^i\to (\zp )^i$ be given by $t_{\bf c}({\bf x})={\bf x}+{\bf c}$. The map  $t_{\bf c}$ induces an isomorphism of $\zp$-modules $( t_{\bf c})_*:\zi \to \zi $

\medskip

\noindent 
{\bf Lemma 6.2.} {\it  Let $\mu  \in  \zi $. Then we have  $( t_{\bf c})_*(\mu )=T_{\bf c}(\mu )$. }

\medskip

\noindent
{\bf Proof.}  Let $\mu =\big( \sum _ {{ k} \in (\zp /p^n\zp )^i}\mu ^{(n)}(k)[k]\big) _{n\in \Nbb}\in \varprojlim  _n \zp [(\Zbb /(p^n))^i]$. Then 
$(t_{\bf c})_*(\mu )$ $=\big( \sum _ {{ k} \in (\zp /p^n\zp )^i}\mu ^{(n)}(k)[k+{\bf  c} ]\big) _{n\in \Nbb}=
\big( \sum _ {{ k} \in (\zp /p^n\zp )^i}\mu ^{(n)}(k - {\bf c})[k]\big) _{n\in \Nbb} =T_{\bf c}(\mu )$.   \hpb

\medskip

Let $d\in \zp ^\times$ and let $\mu \in \zi $. The multiplication by $d$ on $(\zp )^i$ induces an  isomorphism of $\zp$-modules $m_d:\zi \to \zi$.
The proof of the next lemma we leave as an exercise.

\medskip

\noindent
{\bf Lemma 6.3.} {\it Let $d\in \zp ^\times$,  ${\bf c} \in (\zp )^i$ and let $\mu \in \zp [[ (\zp )^i]]$. Then we have
\begin{enumerate}
\item[i)]  $P(m_d(\mu ))(T_1,\ldots ,T_i)=P(\mu )((1+T_1)^d-1, \ldots ,(1+T_i)^d-1)$,

\smallskip

\item[ii)] $  \sum_{x\in (\Zbb /(p^n)^i}     (m_d(\mu ))^{(n)}   (x)[ x]    =\sum _{x\in (\Zbb /(p^n))^i}\mu ^{(n)}(d^{-1}x)[x]\in \zp [ (\Zbb /(p^n))^i]$

\smallskip
\noindent
 for $n\in \Nbb$,

\smallskip

\item[iii)] $\int _{(\zp )^i}x_1^{n_1}x_2^{n_2}\ldots x_i^{n_i}d(m_d(\mu ))(x_1,x_2\ldots x_i)=$
\smallskip

\noindent
$d^{\sum _{k=1}^i  n_k}\int _{(\zp )^i}x_1^{n_1}x_2^{n_2}\ldots x_i^{n_i}d\mu (x_1,x_2\ldots x_i)$,

\smallskip

\item[iv)]  $m_d\circ T_{\bf c} =T_{d\bf c}\circ m_d$.
\hpb
\end{enumerate}     
}

\smallskip

Because of the property ii) of Lemma 6.3  the measure $m _d(\mu )$ we shall denote also by $\mu \circ d^{-1}$. (This notation is also used in \cite[Chapter 2]{L}).

\medskip

Let $x\in \Qbb _p$. We denote by $x\zp$ the $\zp$-submodule of $\Qbb _p$ generated by $x$.
Let $\mu$ be a measure on $(\zp )^i$ with values in $\Qbb _p$. Then there exists $x\in \Qbb _p$ such that the measure $\mu$ has values in  $x\zp$.

 Below we define a kind of exterior product of measures.

\medskip

\noindent 
{\bf Lemma 6.4.} {\it  Let $x_i,y_j,z_k\in \Qbb _p$. Let
\[
\al =\Big( \al ^{(n)}:(\Zbb /(p^n))^i \to x_i\zp \Big)_{n\in \Nbb}\;\;{\text and} \;\; \be =\Big( \be ^{(n)}:(\Zbb /(p^n))^j\to y_j\zp \Big)_{n\in \Nbb}
\]
be measures on $(\zp )^i$ and $(\zp )^j$ respectively. Then
\[
\al \cdot \be :=\Big( (\zpn )^i\times ( \zpn )^j\ni  (a,b)\mapsto   \al ^{(n)}(a)\cdot  \be ^{(n)}(b) \in x_iy_j\zp \Big) _{n\in \Nbb}
\]
is a measure on $(\zp )^{i+j}$. The product of measures $\al \cdot \be$ considered as a measure on $(\zp )^{i+j}$ with values in $\Qbb _p$ does not depend on the choice of $x_i$ and $y_j$.
Moreover if $\ga  $ is a measure on $(\zp )^k$ with values in $z_k\zp$ then
$(\al \cdot \be) \cdot \ga =\al \cdot (\be \cdot \ga)$.  If $\al  \in \zp [[(\zp )^i]]$ and  $\be  \in \zp[[(\zp )^j]]$  then
$P(\al \cdot \be )(T_1,\ldots ,T_{i+j})=P(\al )(T_1,\ldots ,T_i)\cdot P(\be )(T_{i+1},\ldots ,T_{i+j})$.

If $\al$ and $\be$ are distributions with values in $\Qbb _p$ then $\al \cdot \be $ defined as above is also a distribution with values in $\Qbb _p$.}
\hpb

\medskip

If $\al$ is a measure on $(\zp )^i$ then the measure $\al \cdot \al$ on $(\zp  )^{2i}$ we shall also denote by $\al ^2$.  
 
\medskip

\noindent
{\bf Example 6.5. } {\it Let $\si \in G_\Qbb$.  The measure $T_{(\chs ,\chs )} (D_2(\si ))$ is given by
\medskip
\[
 \big( T_{(\chs ,\chs )}(D_2(\si )) \big) ^{(n)}(a,b)=                        \ga _{ \sbn -b}^{(n)}(\si )-\ga _ {\sbn  -a}^{(n)}(\si )+
\]
\[
\left\{
\begin{array}{ll}
\al _{\sbn -a}^{(n)}(\si ) & \text{ if }\;  0<a<b\leq \sbn  \\
                                       & \text{ or }\; \sbn <a<b \leq p^n  \\
                                       & \text{ or }\;  0<b<\sbn <a\leq p^n \\
                                       & \text{ or }\;  b=\sbn <a\leq p^n , \\
-\al _{\sbn -b}^{(n)}(\si ) & \text{ if }\;  0<b<a\leq \sbn\\
                                        & \text{ or }\; \sbn <b<a \leq p^n  \\
                                        & \text{ or }\;  0<a<\sbn <b\leq p^n \\
                                        & \text{ or }\;  a=\sbn <b\leq p^n , \\
0  & \text{ if }\;  a=b
\end{array}
\right.
\]
\medskip
for $(a,b)\in \Zbb /(p^n)\times \Zbb /(p^n)$ .        \hpb
}

\medskip

\noindent
{\bf Example 6.6. } {\it  Let $c=(c_1,c_2)\in (\zp )^2$. Then $\de _{c_1}\cdot \de _{c_2}=\de _c$ in $\zp [[(\zp )^2]]$.}  \hpb

\bigskip

\section{Group theoretical properties of measures}

\medskip

While studying Galois representations on fundamental groups and on torsors of paths one can notice that some results and constructions are purely group theoretical, others are forced by geometrical constrains and still others by number theory (for  example the appearance of the Vandiver conjecture in \cite[Theorem 6]{I}).

In this section we shall work in the group theoretical context. We shall complete Theorem 2.6 from \cite{W11}. We shall also show that the construction of measures in section 3 from projective systems of power series is (almost) tautological.

We recall from section 3 that $V=V_0=\Pbb ^1_{\bQ}\setminus \{ 0,1,\infty \}$ and  $V_n:=\Pbb ^1_ {\overline \Qbb }\setminus (\{ 0,\infty \}\cup \mu _{{p^n}})$.
The covering maps $p_n^{n+m}:V_{n+m}\to V_n,\;  \zfk \mapsto \zfk ^{p^m}$ induce injective morphisms $(p_n^{n+m} )_* :\pi _1(V_{n+m},\01 )\to \pi _1(V_n,\01 )$.  Below we shall identify the group $\pi _1(V_n,\01 )$ with its image by the map $(p_0^n)_*$ in $\pi _1(V ,\01 )$.
Then for each $n$ we have
$$
\pi _1(V_n,\01 )=\ker \big( \pi _1(V ,\01 ) \to \pi _1(\Pbb ^1_{\bQ}\setminus \{ 0,\infty \},\01 )/p^n \pi _1(\Pbb ^1_{\bQ}\setminus \{ 0,\infty \},\01 )\big) ,
$$
where the map $\pi _1(V ,\01 ) \to \pi _1(\Pbb ^1_{\bQ}\setminus \{ 0,\infty \},\01 )$ is induced by the inclusion $V \hookrightarrow \Pbb ^1_{\bQ}\setminus \{ 0,\infty \}$ (see \cite[page 287]{NW}). Hence it follows that
$$
\bigcap _{i=0}^\infty \pi _1 (V_i,\01 )   =\ker \big( \pi _1(V,\01 )\to  \pi _1(\Pbb ^1 _{\overline \Qbb} \setminus \{0,\infty \},\01 )\big) .
$$
Let $g \in \ker \big( \pi _1(V,\01 )\to  \pi _1(\Pbb ^1 _{\overline \Qbb} \setminus \{0,\infty \},\01 )\big) $.
Then for any $n$ we can consider $g$ as an element of $\pi _1(V_n,\01 )$. Hence for any $n$ we get a power series 

\noindent
$E_n (g) (X_n,Y_{0,n},\ldots ,Y_{p^n-1,n})\in  \Qbb  _p\{\{\Yc _n\}\} $ such that $(p_n^{n+m} )_* (E_{n+m}(g))=E_n(g)$ for $n,m\in \Nbb$.

We recall that $\Qbb _p\{\{\Yc _n^0\}\}$ is a $\Qbb _p$-algebra of non-commutative power series on elements of $\Yc _n^0 =\{ Y_{0,n},\ldots ,Y_{p^n-1,n}\}$. We denote by $L(\Yc _n^0)$ the Lie algebra of Lie power series in  $\Qbb _p\{\{\Yc _n^0\}\}$. Let us set
\[
E_n^0 (g) :=E_n (g)(0,Y_{0,n},\ldots ,Y_{p^n-1,n})\; .
\]
Then $E_n^0 (g)\in \Qbb _p\{\{\Yc _n^0\}\}$. The morphism of $\Qbb _p$-algebras $(p_n^{n+m} )_* :  \Qbb  _p\{\{\Yc _{n+m}\}\}\to   \Qbb  _p\{\{\Yc _n\}\}$ induces a morphism of $\Qbb _p$-algebras
$$
\pfk _n^{n+m} :  \Qbb  _p\{\{\Yc _{n+m}^0\}\}\to   \Qbb  _p\{\{\Yc _n^0\}\}
$$
such that $\pfk _n^{n+m}(Y_{i+kp^n,n+m})=Y_{i,n}$ for $0\leq i<p^n$ and $0\leq k<p^m$ because $(p_n^{n+m} )_* (Y_{i+kp^n,n+m})=\exp (-kX_n)\cdot Y_{i,n}\cdot \exp (kX_n)$.

For any $i>0$ and any   $ {\bf a}\in (\Zbb /(p^n))^i$ let us define coefficients    $ K _i ^{(n)}(g)( {\bf a})$ by the identity
$$
 E_n^0 (g)=1+\sum _{i=1}^\infty \Big(\sum _{ {\bf a}\in (\Zbb /(p^n))^i}    K _i ^{(n)}(g)( {\bf a})                {\bf Y}_{\bf a}\Big) ,
$$
where ${\bf a}=(a_1,a_2,\ldots ,a_i)\in  (\Zbb /(p^n))^i$ and $ {\bf Y_a}=Y_{a_1,n}Y_{a_2,n}\ldots Y_{a_i,n}$.

\medskip

\noindent 
{\bf Proposition  7.1.} {\it Let $g \in  \ker \big( \pi _1(V,\01 )\to  \pi _1(\Pbb ^1 _{\overline \Qbb} \setminus \{0,\infty \},\01 )\big) $.
Then we have:
\begin{enumerate}
\item[i)] for each $i>0$ the family of functions
$$
\big(  (\Zbb /(p^n))^i \ni {\bf a} \mapsto   K _i ^{(n)}(g)( {\bf a})    \in \Qbb _p \big) _{n\in \Nbb}
$$ 
is a measure on $(\zp )^i$ with values in ${\frac{1}{i!}}\zp \subset \Qbb _p$,
\item[ii)] for each $n$, the power series $E_n^0(g)$ belongs to the subgroup $\exp (L(\Yc _n^0))$ of $\Qbb_p \{ \{ (\Yc _n^0 \}\}$, hence the coefficients of the power series $E_n^0(g)$ satisfy the shuffle relations.
\end{enumerate}
}

\smallskip

\noindent
{\bf Proof.}
Observe that for any ${\bf a}\in (\zp )^i$ the coefficient of $E_n(g)\in \qp \{\{\Yc _n\}\}$ at ${\bf Y}_{\bf a}$ is equal to the coefficient of  $E_n^0(g)\in \qp \{\{\Yc _n^0\}\}$ at ${\bf Y}_{\bf a}$.
 Hence in order to show that we have a measure it is sufficient  to  repeat the proof of Proposition 2.3  from \cite{W11} in the context of the proposition. 
The worst possible denominator appears in $K _i ^{(n)}(g) (k,k,\ldots ,k)={\frac{1}{i!}}(K  _1^{(n)}(g)(k))^i\in {\frac{1}{i!}}\zp$. Hence the measure 
has values in $ {\frac{1}{i!}}\zp$. 

By the very definition of the embedding $E_n^0$,  we have $E_n^0(y_{i,n})=e^{Y_{i,n}}$. Hence it follows from the Baker-Cambell-Hausdorff theorem that
$E_n^0(g) \in \exp (L(\Yc _n^0))$ (see \cite[Theorem 5.19]{MKS}). Therefore the coefficients of the power series $E_n^0$ satisfy the shuffle relations 
(see \cite[Theorem 2.5.]{R}).   \hpb

\medskip

The measure
$$
\Big( (\Zbb /(p^n))^i \ni {\bf a}\mapsto    K _i ^{(n)}(g)( {\bf a})  \in  {\frac{1}{i!}}\Zbb _p \Big) _{n\in \Nbb }
$$
we denote by ${\bf K} _i(g)$. For $i=0$ we define a measure
$$
{\bf K} _0(g):=\Big( (\Zbb /(p^n))^0 \ni * \mapsto 1 \in \Qbb _p \Big) _{n\in \Nbb }.
$$

\medskip
Let $\si \in  G_\Qbb$.
Then the  measure $\Kc _r(\si )$ from section 3 is the measure ${\bf K} _r (\ffk _\pi (\si ))$ in the notation of this section.

\medskip

If $n=0$ we usually omit the subscript $0$, hence we have an embedding $E:\pi _1(V,\01 )\to \qp \{\{ X,Y\} \}$.

If $\phi \in \Qbb _p\{ \{ X,Y\} \}$
and $w\in \Mc$ is a monomial in non-commuting variable $X$ and $Y$, then 
$\langle \phi ,w\rangle$ is  the coefficient of the series $\phi$ at the monomials $w$.

 \medskip

\noindent 
{\bf Proposition 7.2.} {\it  Let $g \in  \ker \big( \pi _1(V,\01 )\to  \pi _1(\Pbb ^1 _{\overline \Qbb} \setminus \{0,\infty \},\01 )\big) $.  Let us set $w(n_0,n_1,\ldots , n_{r-1},n_r) = X^{n_0}Y X^{n_1}Y\ldots  X^{n_{r-1}}YX^{n_r}$.  The power series $E(g)\in \Qbb _p\{ \{ X,Y\} \}$ and the family of measures $\{ {\bf K} _i(g)\} _{i\in \Nbb }$ determine each other.  
 More precisely we have
\[
\langle E(g), X^{n_0}Y X^{n_1}Y\ldots  X^{n_{r-1}}Y X^{n_r}\rangle =
\]
\[
{\frac{1}{n_0!n_1!\ldots n_r!}}\int _{(\Zbb _p)^r}(-x_1)^{n_0}(x_1-x_2)^{n_1}\ldots (x_{r-1}-x_r)^{n_{r-1}}x_r^{n_r}d{\bf K} _r(g)(x_1,x_2,\ldots ,x_r)
\]
and
\[
F(           { \bf K} _r(g))(X_1,X_2,\ldots ,X_r)=
\]
\[
\sum _{{ \mathbf n}\in \Nbb ^r}\langle E(g),w(n_0,..n_{r-1},0) \rangle (-\sum _{i=1}^rX_i)^{n_0} (-\sum _{i=2}^rX_i)^{n_1}\ldots (-X_{r-1}-X_r)^{n_{r-2}}(-X_r)^{n_{r-1}} ,
\]
where ${\mathbf n} =(n_0,n_1,\ldots  ,n_{r-1})$.
}

\smallskip

\noindent
{\bf Proof.} The first formula is proved in the same way  as \cite[Theorem 2.6.]{W11}.

In order to show the second formula let us notice that by the very definition we have
$F(  {\bf K} _r(g)   )(X_1,\ldots ,X_r)=$
$$\sum _{  {\mathbf n}   \in (\Nbb )^r         }   {\frac{1}{n_0!\ldots n_{r-1}!}}
 \Big( \int _{(\Zbb _p)^r }x_1^{n_0}\ldots x_r^{n_{r-1}}d  {\bf K} _r(g) (x_1,\ldots ,x_r)\Big) X_1^{n_0}\ldots  X_r^{n_{r-1}}\, .
$$
Observe that we have the following equalities of formal power series in $\qp [[ X_1,\ldots ,X_r]]$
$$
\sum _{   {\mathbf n}   \in (\Nbb )^r            }   {\frac{1}{n_0!\ldots n_{r-1}!}}(-x_1)^{n_0} (x_1 -x_2)^{n_1}\ldots (x_{r-1} -x_r)^{n_{r-1}}
\big(-\sum _{i=1}^rX_i\big)^{n_0}\big(-\sum _{i=2}^rX_i\big)^{n_1}\ldots \big( -X_r\big) ^{n_{r-1}}=
$$

$$
e^{x_1( \sum _{i=1}^rX_i)}e^{(x_2-x_1)(\sum _{i=2}^rX_i)}\ldots e^{(x_r-x_{r-1})X_r}=e^{x_1X_1+x_2X_2+\ldots  +x_rX_r}=
$$

$$
\sum _{n_0=0}^\infty \ldots \sum _{n_{r-1}=0}^\infty  {\frac{1}{n_0!\ldots n_{r-1}!}}x_1^{n_0}\ldots x_r^{n_{r-1}}  X_1^{n_0}\ldots  X_r^ {n_{r-1}}\, .
$$
Integrating over $(\Zbb _p)^r$ the coefficients of the two power series we get the second formula.
\hpb

\medskip

We denote by
$
{\rm Dist} \big( (\zp )^m,\Qbb _p\big)
$
the $\Qbb _p$-module of distributions on $(\zp )^m$ with values in $\Qbb _p$. Let $x\in \Qbb _p$. We denote by 
$
{\rm Meas}\big(  (\zp )^m,\Qbb _p\big) \;\;\; ({\rm resp.}\; {\rm Meas}\big(  (\zp )^m,x\zp \big)
$
the $\Qbb _p$-module (resp. the $\zp$-module) of measures on $(\zp )^m$ with values in $\Qbb _p$ (resp. in $x\zp$).
The $\Qbb _p$-module $ \prod _{m=0}^\infty {\rm Dist} \big( (\zp )^m,\Qbb _p\big)$ we equipped with a product $*$ defined by
 $(\al _i)^\infty _{i=0} * (\be _i)^\infty _{i=0}:=(c _i)^\infty _{i=0}$, where $c_i:=\sum _{j+k=i}\al _j\cdot \be _k$
and where $\al _j\cdot \be _k$ is the product of distributions introduced in Lemma 6.4.
Then the product  $ \prod _{m=0}^\infty {\rm Dist} \big( (\zp )^m,\Qbb _p\big)$ is also a $\Qbb _p$-algebra.
Observe also  that  $\prod _{r=0}^\infty \qp$ acts on  $ \prod _{m=0}^\infty {\rm Dist} \big( (\zp )^m,\Qbb _p\big)$ component-wisely.

Let $f\in \qp \{ \{\Yc _n^0\} \}$. We denote by $h_r(f)$ the homogenous term of degree $r$ of $f$. Then $f=\sum _{r=0}^\infty h_r(f)$. The term $h_0(f)$ is the constant term of $f$. We define an action of $\prod _{r=0}^\infty \qp$ on $ \{ \{\Yc _n^0\} \}$. 
 Let $c =(c_r)_{r\in\Nbb}\in \prod _{r=0}^\infty \qp$.
 Then we set $c(f):= \sum _{r=0}^\infty c_r h_r(f)$.

Let $\varprojlim  _n \Qbb _p\{\{\Yc _n^0\}\}$ be the inverse limit taken with respect to maps 

\noindent
$\pfk _n^m:\Qbb _p\{\{\Yc _m^0\}\}\to \Qbb _p\{\{\Yc _n^0\}\}$ for $m\geq n$.
The inverse limit is equipped with the product induced by the multiplication of series in each $ \Qbb _p\{\{\Yc _n^0\}\}$.
Hence $\varprojlim  _n \Qbb _p\{\{\Yc _n^0\}\}$  is a $\Qbb _p$-algebra.
 Let $c  \in \prod _{r=0}^\infty \qp$.
 If $(f_n)_{n\in \Nbb}\in \varprojlim  _n \Qbb _p\{\{\Yc _n^0\}\}$ then
$(c(f_n))_{n\in \Nbb}\in \varprojlim  _n \Qbb _p\{\{\Yc _n^0\}\}$.

\medskip

\noindent 
{\bf Proposition 7.3.} {\it  The natural map
$$
 \varprojlim  _n \Qbb _p\{\{\Yc _n^0\}\} \to \prod _{m=0}^\infty  {\rm Dist} \big( (\zp )^m,\Qbb _p\big)
$$
is an isomorphism of $\Qbb _p$-algebras.  Moreover the above isomorphism is compatible with the actions of $ \prod _{r=0}^\infty \qp$ on both $\qp$-algebras.      } 

\smallskip

\noindent
{\bf Proof.}  It follows from \cite[E III.5 Corollaire 2]{B} that we have a bijection
\[
 \varprojlim  _n \big( \prod _{r=0}^\infty  \qp [(\Zbb /(p^n))^r]\big) \cong  \prod _{r=0}^\infty \big(  \varprojlim  _n  \qp [(\Zbb /(p^n))^r] \big ) ,
\]
where the maps of projective systems are induced by the reduction maps from $(\Zbb /(p^n))^{r+1}$ to $(\Zbb /(p^n))^r$. The bijection is easily seen to be an isomorphism of $\qp$-modules. The group algebra we identify with a $\qp$-module of homogenous polynomials of degree $r$ in non-commuting variables $Y_{0,n},\ldots ,Y_{{p^n-1},n}$ by the correspondence

$$ 
  (\Zbb /(p^n))^r \ni (i_1,\ldots ,i_r)\mapsto  Y_{i_1,n},\ldots ,      Y_{i_r,n}         \, .
$$

\smallskip

\noindent
Then $\prod _{r=0}^\infty  \qp [(\Zbb /(p^n))^r]$ is identified with the $\qp$-module $\Qbb _p\{\{\Yc _n^0\}\}$.
On the other side
$$
\varprojlim  _n  \qp [(\Zbb /(p^n))^r] =  {\rm Dist} \big( (\zp )^r,\Qbb _p\big)
$$
by the very definition. Hence we get an isomorphism of $\qp$-modules
$$
 \varprojlim  _n \Qbb _p\{\{\Yc _n^0\}\} \cong \prod _{r=0}^ \infty  {\rm Dist} \big( (\zp )^r,\Qbb _p\big)\, ,
$$
which is also compatible with the action of $\prod  _{r=0}^\infty    \qp$ on both sides, and which we denote by $\Ic$.

Let $f=(f_n)_{n\in \Nbb}$, $g=(g_n)_{n\in \Nbb}$ $\in  \varprojlim  _n \Qbb _p\{\{\Yc _n^0\}\} $. Then $f_n\cdot g_n =$

\noindent
$\sum _{r=0}^\infty \big( \sum _{i+j=r}h_i(f_n)\cdot h_j(g_n)\big)$ and $\sum _{i+j=r}h_i(f_n)\cdot h_j(g_n)$ correspond to the $r$-th component of $\Ic (f) *  \Ic (g)$.
\hpb

\medskip

 We define a map
\[
{\bf K }:\ker \big( \pi _1(V,\01 )\to  \pi _1(\Pbb ^1 _{\overline \Qbb} \setminus \{0,\infty \},\01 )\big) \to \prod _{m=0}^\infty   {\rm Meas} \big( (\zp )^m,{\frac{1}{m!}}\zp \big)
\]
by ${\bf K}(g):=({\bf K} _m(g))^\infty _{m=0}$.

\medskip

\noindent 
{\bf Proposition 7.4.} {\it 
\begin{enumerate}
\item[i)]  The map ${\bf K}$ is injective.

\item[(ii)]  Let $   g,h \in \ker \big( \pi _1(V,\01 )\to  \pi _1(\Pbb ^1 _{\overline \Qbb} \setminus \{0,\infty \},\01 )\big) $. Then ${  {\bf K}  }(gh)= { {\bf K} }(g)* { {\bf K} }(h)$.
\end{enumerate}
}

\noindent
{\bf Proof.}  The map $E: \pi _1(V,\01 )\to \Qbb _p\{\{ X,Y\}\}$ is injective.  Hence it follows from Proposition 7.2 that the map $\be$ is injective.

Let $ g,h \in \ker \big( \pi _1(V,\01 )\to  \pi _1(\Pbb ^1 _{\overline \Qbb} \setminus \{0,\infty \},\01 )\big) $. Then for each $n$ we have
$E_n^0(g  h) =E_n^0(g) E_n^0(h)$. Hence it follows from Proposition 7.3 that ${  {\bf K} }(gh)= {\bf K} (g)*  {\bf K}  (h)$. \hpb

\medskip

\noindent 
{\bf Question  7.5.} What is the image of the map $\bf K$?

\medskip

At the end of this section we return to the Galois context.
We recall that $\pi$ is a path on $V$ from $\01$ to $\10$ (the real interval $[0,1]$). The map
\[
\ffk _\pi :G_{\Qbb}\to  \pi _1(V,\01 ),\; \si \mapsto \ffk _\pi (\si )=\pi ^{-1}\cdot \si \cdot \pi \cdot \si ^{-1}
\]
or equivalently the map $G_\Qbb \ni \si \mapsto  E(\ffk _\pi (\si ))\in \Qbb _p\{\{ X,Y\}\}$ 
is a cocycle. Its image satisfy the Deligne-Drinfeld-Ihara relations (i.e, $\Zbb /(2)$, $\Zbb /(3)$ and  $\Zbb /(5)$ symmetries) in the pro-$p$ group  $\pi _1(V,\01 )$ or in the $\Qbb _p$-algebra $ \Qbb _p\{\{ X,Y\}\}$  .

\medskip

\noindent 
{\bf Question  7.6.} What are the  relations which  satisfies the image of the composition ${\bf K} \circ \ffk _\pi$ in $  \prod _{m=0}^\infty   {\rm Meas} \big( (\zp )^m,{\frac{1}{m!}}\zp \big) $ ?

\medskip

For example  in degree $1$ we have
\[
{\bf K}    _1(\ffk _\pi (\si )) - { \bf K} _1(\ffk _\pi (\si )) \circ (-1) =E_{1,\chi (\si )}+{\frac{1-\chi (\si )}{2}}\de _0 \, .
\]
(see \cite[Proposition 7]{W12} and also
 \cite[Proposition 5.13]{NW2}).

\bigskip

\section{Symmetries of the measures  on $\zp$ and $(\zp )^2$ }

\medskip
In this section we shall study only the measures $\Kc _1(\si )$ and $\Kc _2(\si )$ for $\si \in G _ \Qbb$.  We use notations from the very beginning of section 4.
For more clarity we denote the measures $\Kc _ 1(\si )$ and $\Kc _2(\si )$ by $\al (\si )$ and $\be  (\si )$ respectively, i.e.
$$
\al (\si ):=\big(\sum _{i=0}^{p^n-1}\al _i^{(n)}(\si )[i]\big) _{n\in \Nbb}\in  \varprojlim _n \zp [\Zbb /(p^n)]
$$
and 
$$
\be (\si ):=\big( \sum _{0\leq a,b <p^n}\be _{a,b}^{(n)}(\si )[a,b]\big) _{n\in \Nbb}\in  \varprojlim _n \zp [ ( \Zbb /(p^n) )^2 ] \, .
$$

\noindent 
{\bf Proposition  8.1.} {\it Let $\si  \in G_\Qbb$. Then we have
\[
\al (\si ) -\al (\si )\circ (-1)+ T_{\chs } (\al (\si )\circ (-1) )    -   T_{\chs }(\al (\si )) +
\]
\[
 \Mc (\chs ) +{\frac{1-\chs }{2}}\de _0 +{\frac{1-\chs }{2}}\de _{\chs }=0\,.
\]
Let $\si  \in G_{\Qbb (\mu _\infty )}$. Then we have
\[
\al (\si ) -\al (\si )\circ (-1)+ T_{1 } (\al (\si )\circ (-1) )    -   T_{1 }(\al (\si )) =0\, .
\]
}

\noindent
{\bf Proof.}  The first formula of the  proposition follows immediately from Corollary 4.4, from the definition of the measure $ \Mc (\chs ) $ in section 5 and from the definition of operations on measures in section 6. The second formula follows from the fact that $\chs =1$ for $\si \in   G_{ \Qbb  (\mu _{p^\infty})}$.     \hpb

\medskip

\noindent 
{\bf Theorem  8.2.} {\it Let $\si \in G_{ \Qbb  (\mu _{p^\infty})}$. Let $(a,b)\in \Zbb /p^n \times \Zbb /p^n$.  Then we have 
\begin{enumerate}
\item[i)]  

\noindent
$\be _{a,b}^{(n)}(\si ) -  \be _{-a,-b}^{(n)}(\si ) + \be _{1 -a,1 -b}^{(n)}(\si ) - \be _{a-1,b-1}^{(n)}(\si ) +D_2^{(n)}(\si ) (a,b) -
$

\noindent
$
D_2^{(n)}(\si ) (a,b) - D_2^{(n)} (\si ) (a-1,b-1) +\al _{1-a}^{(n)}(\si )\cdot \de _1^{(n)}(b) -  \de _1^{(n)}(b)\cdot \al _{1-b}^{(n)}(\si )=0$,

\smallskip

\item[ii)] $\be (\si ) -\be (\si )\circ (-1) +T_{(1,1)}(\be (\si )\circ (-1) )-          T_{(1,1)}(\be (\si ) ) +$

\noindent
$D_2(\si ) - T_{(1,1)}(D_2(\si ))
+ \big( T_{1 } (\al (\si )\circ (-1) ) \big) \cdot \de _1- \de _1\cdot  \big( T_{1 } (\al (\si )\circ (-1) ) \big)=0$.

\end{enumerate}
}

\medskip

\noindent
{\bf Proof.} Let us consider terms $A_n$, $B_n$,\ldots , $K_n$ from section 4.
 Observe that $B_n=D_n=G_n=J_n=K_n=1$ for  $\si \in G_{ \Qbb  (\mu _{p^\infty})}$ because then  $\chs =1$.  Hence it follows from Proposition 4.3 that for  $\si \in G_{ \Qbb  (\mu _{p^\infty})}$ we have  $H_n \cdot F_n \cdot C_n \cdot  A_n \equiv 1$  modulo $I_n ^3$. 

We shall write terms in degree 2 in variables $Y_i$'s of the product $H_n \cdot F_n \cdot C_n \cdot  A_n $ starting from $A_n$ to $H_n$ and next from the products $C_nA_n$, $F_nA_n$, $ H_nA_n$, $ F_nC_n$, $H_nC_n$ and $H_nF_n$.
It follows from \cite[Lemma 4.1]{W11} that for $\si \in  G_{ \Qbb  (\mu _{p^\infty})}$, $\al _i^{(n)}(\si )=\al _{p^n - i}^{(n)}(\si )$. One then sees that the terms of the form $\al _a^{(n)}(\si )\al _b^{(n)}(\si )Y_aY_b$ disappear. The remaining terms are
$$
\sa  \be _{a,b}^{(n)}   (\si )Y_aY_b 
 -    \sum _{b=1}^{p^n-1} \al _{p^n-b} ^{(n)} (\si )  \big( \sum _{b<a\leq p^n}Y_a\big) Y_b+ \sum _{a=1}^{p^n-1} \al _{p^n-a} ^{(n)} (\si )Y_a \big( \sum _{a<b\leq p^n}Y_b\big) +
$$
$$
(-1)\sa \be _{-a,-b }^{(n)} (\si )Y_aY_b +
\sbp \ga _{ -b }^{(n)} (\si ) \big(   \sap Y_a\big) Y_b-  \sap \ga _{ -a }^{(n)} (\si )Y_a \big(  \sbp Y_b\big)  +
$$
$$
 \sum _{ a= 2}^{p^n-1}\al _{1 -a}^{(n)}(\si )\Big(  - Y_a \big(\sum _{j=a+1}^{p^n}Y_j\big)+  \big(\sum _{j=a+1}^{p^n}Y_j\big) Y_a\Big) + 
 \sa \be _{1 -a,1 -b}^{(n)}(\si )Y_aY_b +
$$
$$
  ( -1) \sum _{j=0}^{p^n-1} \ga _{1 -j}^{(n)}(\si )\Big( \big(\sum _{k=0}^{p^n-1}Y_k\big)Y_j - Y_j   \big(\sum _{k=0}^{p^n-1}Y_k\big) \Big)\, 
-\sa \be _{a-1 ,b- 1 }^{(n)}(\si )Y_aY_b \,.
$$

\medskip

The terms 
$$
\sa  \be _{a,b}^{(n)}   (\si )Y_aY_b  - \sa \be _{-a,-b }^{(n)} (\si )Y_aY_b +  \sa \be _{1 -a,1 -b}^{(n)}(\si )Y_aY_b 
-\sa \be _{a-1 ,b- 1 }^{(n)}(\si )Y_aY_b 
$$
form the measure
$$
\be _2(\si ) -\be _2(\si )\circ (-1) +T_{(1,1)}(\be _2(\si )\circ (-1) )-          T_{(1,1)}(\be _2(\si ) ) \, .
$$
Observe that 
$$
\sbp \ga _{ -b }^{(n)} (\si ) \big(   \sap Y_a\big) Y_b-  \sap \ga _{ -a }^{(n)} (\si )Y_a \big(  \sbp Y_b\big)  +
\sum _{a=1}^{p^n-1} \al _{p^n-a} ^{(n)} (\si )Y_a \big( \sum _{a<b\leq p^n}Y_b\big) +
$$
$$
( - 1)   \sum _{b=1}^{p^n-1} \al _{p^n-b} ^{(n)} (\si )  \big( \sum _{b<a\leq p^n}Y_a\big) Y_b=\sa  D_2^{(n)}([a,b])Y_aY_b\, .
$$
The remaining  terms are
$$
 - \sum _{j=0}^{p^n-1} \ga _{1 -j}^{(n)}(\si )\Big( \big(\sum _{k=0}^{p^n-1}Y_k\big)Y_j - Y_j   \big(\sum _{k=0}^{p^n-1}Y_k\big) \Big) +
$$
$$
\sum _{ a= 2}^{p^n-1}\al _{1 -a}^{(n)}(\si )\Big(  - Y_a \big(\sum _{j=a+1}^{p^n}Y_j\big)+  \big(\sum _{j=a+1}^{p^n}Y_j\big) Y_a\Big) =
$$ 
$$
\sa \big( - \ga _{1 -b}^{(n)}(\si ) +  \ga _{1 -a}^{(n)}(\si )\big) Y_aY_b-\sum _{1<a<b\leq p^n}\al _{1-a}^{(n)}(\si )Y_aY_b +
 \sum _{1<b<a\leq p^n}\al _{1-b}^{(n)}(\si )Y_aY_b +
$$
$$
(-1) \sum _{1<a \leq p^n}\al _{1-a}^{(n)}(\si )Y_aY_1 + \sum _{1<b \leq p^n}\al _{1-b}^{(n)}(\si )Y_1Y_b +
\sum _{1<a \leq p^n}\al _{1-a}^{(n)}(\si )Y_aY_1 - \sum _{1<b \leq p^n}\al _{1-b}^{(n)}(\si )Y_1Y_b \, .
$$
The last two sums form a measure $T_1(\al (\si )\circ (-1)) \cdot \de _1 - \de _1 \cdot T_1(\al (\si )\circ (-1))$, while the rest form a measure $-T_{(1,1)}(D_2(\si ))$ (see Example 6.5 for  $\si \in G_{ \Qbb  (\mu _{p^\infty})}$). \hpb

\medskip

For   $\si \in G_{  \Qbb}$ the corresponding formula is much more complicated.  First we recall the relations between measures which will be useful to simplify the calculations.

\smallskip

\noindent 
{\bf Lemma 8.3. }{\it        Let $c\in \zp ^\times$ and let $\si \in G_{\Qbb}$.  There are the following relations between measures     
\begin{enumerate}
\item[i)] $E_{1,c}+E_{1,c}\circ (-1)=(c-1)\delta _0$,
\item[ii)] $ T_c(E_{1,c}  )=E_{1,c} +\Mc (c) +(1-c)\de _0$,
\item[iii)] $\al (\si )- \al (\si )\circ (-1)=E_{1,\chs }+{\frac{1-\chs }{2}}\de _0$,
\item[iv)]  $ T_\chs (\al (\si )  ) -T_{\chs }(\al (\si )\circ(-1))=E_{1,\chs} +\Mc (\chs) +(1-\chs)\de _0 +{\frac{1-\chs }{2}}\de _{\chs } $.
 \end{enumerate}
  }

\medskip

\noindent
{\bf Proof.}  The point i) is  instantaneous. The point iii) follows from \cite[Proposition 7]{W12}. To show  the point ii) we compare the power series corresponding to measures  $T_c(E_{1,c}  )$ and $E_{1,c} +\Mc (c) +(1-c)\de _0$. The point iv) follows from ii) and iii).  \hpb

\medskip

\noindent 
{\bf Theorem  8.4.} {\it   Let $\si \in G_{ \Qbb}$.    Then   in $\zp [[(\zp )^2]]$ we have the following identity 

\smallskip
 \[
\be (\si )- \be (\si )\circ (-1) +T_{( \chs , \chs )} (\be (\si )\circ (-1)) - T_{(\chs ,\chs )}( \be (\si )) +
\]
\[
(-1) \big( \al (\si )\circ (-1) \big) \cdot  E_{1,\chs}- \big(  \al (\si )\circ (-1) \big) \cdot (1-\chs)\de _0  -    E_{1,\chs}   \cdot E_{1,\chs}     + 
\]
\[
 (  -1)       (1-\chs)\de _0  \cdot E_{1,\chs}-    E_{1,\chs}    \cdot          (1-\chs)\de _0 + \big(  T_{\chs } (\al (\si )  ) \big)    \cdot E_{1,\chs}+
\]
\[
 \big( T_{\chs} (\al (\si )  ) \big) \cdot \Mc (\chs ) +                    \big( T_{\chs} (\al (\si )  ) \big) \cdot  (1-\chs)\de _0 - {\frac{7}{8}}  (1-\chs)^2 \de _0 \cdot     \de _0 +D_2 (\si ) +
\]
\[
{\frac{ (1-\chs)^2}{8}}\de _{\chs } \cdot   \de _{\chs }+ {\frac{1-\chs }{2}}\de _0 \cdot \al (\si ) + {\frac{1-\chs }{2}}\de _{\chs}  \cdot \big( T_{\chs }(\al (\si )\circ(-1))\big) +
\]
\[
( -1)  T_{(\chs ,\chs ) }(D_2(\si ))- {\frac{1}{2}}\Mc (\chs)\cdot \Mc (\chs) +\big( T_{\chs }(\al (\si )\circ(-1)) \big) \cdot     \de _{\chs} +
\]
\[
   ( -1)  \de _{\chs}  \cdot                             \big(         T_{\chs }(\al (\si )\circ(-1))    \big)              -  E_{1,\chs}\cdot \Mc (\chs) -       (1-\chs)\de _0 \cdot        \Mc (\chs)  +\frac{1}{2} \Nc _2 (\chs )   =0 \, .
\]
}

\medskip

\noindent
{\bf Proof.}  Let us consider terms $A_n$, $B_n$, ..., $K_n$ from section 4.
We shall write terms in degree 2 in variables $Y_i$'s of the product $K_n \cdot J_n\cdot H_n\cdot G_n\cdot F_n\cdot D_n \cdot C_n \cdot B_n\cdot  A_n $ starting from $A_n$ to $K_n$ and next from the products $B_nA_n$,  $C_nA_n$ to $K_nH_n$.
We then get
\smallskip

$$
 \sa  \be _{a,b}^{(n)}   (\si )Y_aY_b + 
 {\frac{(1-\chi (\si ) )^2}{8}}Y_0 ^2 +
$$
$$ 
  (-1)    \sum _{b=1}^{p^n-1} \al _{p^n-b} ^{(n)} (\si )  \big( \sum _{b<a\leq p^n}Y_a\big) Y_b+ \sum _{a=1}^{p^n-1} \al _{p^n-a} ^{(n)} (\si )Y_a \big( \sum _{a<b\leq p^n}Y_b\big) -\sa \be _{-a,-b }^{(n)} (\si )Y_aY_b +
$$
\[
\sbp \ga _{ -b }^{(n)} (\si ) \big(  \sap Y_a\big) Y_b-  \sap \ga _{ -a }^{(n)} (\si )Y_a \big(  \sbp Y_b\big) +  \big( \sip \al _{p^n-i} ^{(n)} (\si )Y_ i \big) ^2 \, +
\]
\[
  {\frac{\chs - \ls}{2  p^n}}  \big( \sum _{1\leq b<a\leq p^n}Y_aY_b \big) +
  ( -1)  {\frac{\chs - \ls}{2  p^n}} \big( \sum _{1\leq a<b\leq p^n}Y_aY_b \big) + 
\]

$$
 {\frac{1}{2}}\big(  {\frac{\chs - \ls}{  p^n}}\big) ^2 \big( \sa Y_aY_b \big)+
\sum _{a=1}^{\ls}\al _{\ls -a}^{(n)}(\si )   Y_a \big( \sum _{b=1}^{a-1}Y_b\big) -   \sum _{b=1}^{\ls}\al _{\ls -b}^{(n)}(\si )            \big( \sum _{a=1}^{b-1}Y_a\big) Y_b +  
$$

\[
 \sum _{ a=\ls +   1}^{p^n-1}\al _{\ls -a}^{(n)}(\si )\Big(  - Y_a \big( \sum _{j=a+1}^{p^n}Y_j\big)+  \big(\sum _{j=a+1}^{p^n}Y_j\big) Y_a\Big) + 
\]
\[
 \sa \be _{\ls -a,\ls -b}^{(n)}(\si )Y_aY_b - \sum _{j=0}^{p^n-1} \ga _{\ls -j}^{(n)}(\si )\Big( \big(\sum _{k=0}^{p^n-1}Y_k\big)Y_j - Y_j   \big(\sum _{k=0}^{p^n-1}Y_k\big) \Big)\, +
\]
 
\[
       {\frac{1}{2}}  \big(   {\frac{1-\chi (\si )}{2}} \big) ^2(Y_{\ls })^2+ {\frac{1-\chi (\si )}{2}}\Big( Y_{\ls }( \sum _{i=1}^{\ls - 1}Y_i) - ( \sum _{i=1}^{\ls - 1}Y_i) Y_{\ls }\Big)\,+
\]
$$
(-1)   \sip \al _{i- \ls }^{(n)}(\si )\Big(  Y_i ( \sum _{j=1}^{\ls -1}Y_j )-  ( \sum _{j=1}^{\ls -1}Y_j )Y_i \Big)
-\sa \be _{a-\ls ,b- \ls }^{(n)}(\si )Y_aY_b  +
$$
$$
  \big( \sip  \al _{i- \ls }^{(n)}(\si ) Y_i \big) ^2 \, 
 +  {\frac{1}{2}} \sum _{1\leq j < i< \ls }\big(  Y_iY_j - Y_jY_i \big)+ {\frac{1}{2}} \Big(  \sum _{i=1}^{\ls - 1}Y_i  \Big) ^2\,  +
$$

\medskip

$$
{\frac{1-\chi (\si )}{2}}Y_0 \cdot         ( \sum _{b=0}^{p^n-1} \al _b ^{(n)} (\si )Y_ b)
- (  \sum _{a=0}^{p^n-1} \al _{p^n-a} ^{(n)} (\si )Y_ a  )\cdot         (\sum _{b=0}^{p^n-1} \al _b^{(n)} (\si )Y_ b)+
$$
$$
  (  {\frac{\chs - \ls}{ p^n}}(\sap Y_a) )\cdot  ( \sum _{b=0}^{p^n-1} \al _b ^{(n)} (\si )Y_ b)+
(  \sum _{a=0}^{p^n -1}\al _{\ls -a}^{(n)}(\si ) Y_a     )\cdot          (\sum _{b=0}^{p^n-1} \al _b ^{(n)} (\si )Y_ b)+
$$
$$
 ({\frac{1-\chi (\si )}{2}}Y_{\ls }  )\cdot              (\sum _{b=0}^{p^n-1} \al _b ^{(n)} (\si )Y_ b) -
  ( \sap \al _{a- \ls }^{(n)}(\si )  Y_a )\cdot        ( \sum _{b=0}^{p^n-1} \al _b ^{(n)} (\si )Y_ b)+
$$
$$
 (\sum _{a=1}^{\ls -1}Y_a)\cdot       ( \sum _{b=0}^{p^n-1} \al _b^{(n)} (\si )Y_ b) 
- (  \sum _{a=0}^{p^n-1} \al _{p^n-a} ^{(n)} (\si )Y_ a )\cdot    {\frac{1-\chi (\si )}{2}}Y_0      +
$$
$$
  (  {\frac{\chs - \ls}{ p^n}}(\sap Y_a) )\cdot        {\frac{1-\chi (\si )}{2}}Y_0     +
(  \sum _{a=0}^{p^n -1}\al _{\ls -a}^{(n)}(\si ) Y_a     )\cdot    {\frac{1-\chi (\si )}{2}}Y_0        +
$$
$$
 ({\frac{1-\chi (\si )}{2}}Y_{\ls }  )\cdot             {\frac{1-\chi (\si )}{2}}Y_0  -
  ( \sap \al _{a- \ls }^{(n)}(\si )  Y_a )\cdot      {\frac{1-\chi (\si )}{2}}Y_0    +
$$
$$
 (\sum _{a=1}^{\ls -1}Y_a)\cdot      {\frac{1-\chi (\si )}{2}}Y_0 +
 $$
 
$$
 (- 1) (  {\frac{\chs - \ls}{ p^n}}(\sap Y_a) )\cdot       (  \sum _{b=0}^{p^n-1} \al _{p^n-b} ^{(n)} (\si )Y_ b  )        
    - (  \sum _{a=0}^{p^n -1}\al _{\ls -a}^{(n)}(\si ) Y_a     )\cdot        (  \sum _{b=0}^{p^n-1} \al _{p^n-b} ^{(n)} (\si )Y_ b  )   +
$$
$$
(-1)  ({\frac{1-\chi (\si )}{2}}Y_{\ls }  )\cdot          (  \sum _{b=0}^{p^n-1} \al _{p^n-b} ^{(n)} (\si )Y_ b  )      +
  ( \sap \al _{a- \ls }^{(n)}(\si )  Y_a )\cdot       (  \sum _{b=0}^{p^n-1} \al _{p^n-b} ^{(n)} (\si )Y_ b  )             +
$$
$$
( -1)  (\sum _{a=1}^{\ls -1}Y_a)\cdot      (  \sum _{b=0}^{p^n-1} \al _{p^n-b} ^{(n)} (\si )Y_ b  ) +
 (  \sum _{a=0}^{p^n -1}\al _{\ls -a}^{(n)}(\si ) Y_a     )\cdot        {\frac{\chs - \ls}{ p^n}}(\sbp Y_b)       +
$$
$$
 ({\frac{1-\chi (\si )}{2}}Y_{\ls }  )\cdot        {\frac{\chs - \ls}{ p^n}}(\sbp Y_b)        -
  ( \sap \al _{a- \ls }^{(n)}(\si )  Y_a )\cdot        {\frac{\chs - \ls}{ p^n}}(\sbp Y_b)             +
$$
$$
 (\sum _{a=1}^{\ls -1}Y_a)\cdot       {\frac{\chs - \ls}{ p^n}}(\sbp Y_b)   +
 $$
 
$$
 ({\frac{1-\chi (\si )}{2}}Y_{\ls }  )\cdot        (  \sum _{b=0}^{p^n -1}\al _{\ls -b}^{(n)}(\si ) Y_b     )      - 
  ( \sap \al _{a- \ls }^{(n)}(\si )  Y_a)\cdot        (  \sum _{b=0}^{p^n -1}\al _{\ls -b}^{(n)}(\si ) Y_b     )             +
$$
$$
 (\sum _{a=1}^{\ls -1}Y_a)\cdot      (  \sum _{b=0}^{p^n -1}\al _{\ls -b}^{(n)}(\si ) Y_b     )  +
 $$

$$
( -1) ( \sap \al _{a- \ls }^{(n)}(\si )  Y_a )\cdot        ({\frac{1-\chi (\si )}{2}}Y_{\ls }  )        +
 (\sum _{a=1}^{\ls -1}Y_a)\cdot     ({\frac{1-\chi (\si )}{2}}Y_{\ls }  ) +
 $$

$$
(-1)  (\sum _{a=1}^{\ls -1}Y_a)\cdot                                     ( \sbp \al _{b- \ls }^{(n)}(\si )  Y_b ) =0
$$
by Lemma 4.3.
The terms

$$
 \sa  \be _{a,b}^{(n)}   (\si )Y_aY_b  -\sa \be _{-a,-b }^{(n)} (\si )Y_aY_b + \sa \be _{\ls -a,\ls -b}^{(n)}(\si )Y_aY_b +
$$
$$
(-1)\sa \be _{a-\ls ,b- \ls }^{(n)}(\si )Y_aY_b 
$$
form a measure

\begin{equation}
\be (\si )- \be (\si )\circ (-1) +T_{( \chs ,\chs )} (\be (\si )\circ (-1)) - T_{(\chs ,\chs ) }( \be (\si )) .
\tag*{(8.4.1)}
\end{equation} 

The terms
$$
\big( \sip \al _{p^n-i} ^{(n)} (\si )Y_ i \big) ^2 \, + 
 \big( \sip  \al _{i- \ls }^{(n)}(\si ) Y_i \big) ^2 - (  \sum _{a=0}^{p^n-1} \al _{p^n-a} ^{(n)} (\si )Y_ a  )\cdot         (\sum _{b=0}^{p^n-1} \al _b^{(n)} (\si )Y_ b)+
$$
$$
(  \sum _{a=0}^{p^n -1}\al _{\ls -a}^{(n)}(\si ) Y_a     )\cdot          (\sum _{b=0}^{p^n-1} \al _b ^{(n)} (\si )Y_ b) -
  ( \sap \al _{a- \ls }^{(n)}(\si )  Y_a )\cdot        ( \sum _{b=0}^{p^n-1} \al _b ^{(n)} (\si )Y_ b) +
$$
$$
(-1) (  \sum _{a=0}^{p^n -1}\al _{\ls -a}^{(n)}(\si ) Y_a     )\cdot        (  \sum _{b=0}^{p^n-1} \al _{p^n-b} ^{(n)} (\si )Y_ b  )  +  ( \sap \al _{a- \ls }^{(n)}(\si )  Y_a )\cdot       (  \sum _{b=0}^{p^n-1} \al _{p^n-b} ^{(n)} (\si )Y_ b  )        +   
$$
$$
 ( -1)  ( \sap \al _{a- \ls }^{(n)}(\si )  Y_a)\cdot        (  \sum _{b=0}^{p^n -1}\al _{\ls -b}^{(n)}(\si ) Y_b     )   
$$

\noindent
form a measure 

\begin{equation}
( \al  \circ (-1))  \cdot ( \al \circ (-1))+ 
T_\chs ( \al )\cdot T_\chs ( \al )  -  ( \al \circ (-1))\cdot \al +
\tag*{(8.4.2)}
\end{equation}
$$   ( T_\chs ( \al \circ (-1)))\cdot \al -  T_\chs ( \al )\cdot \al  -  ( T_\chs ( \al \circ (-1)))\cdot (\al \circ (-1)) + $$
$$ T_\chs ( \al )\cdot   ( \al \circ (-1)) -  T_\chs ( \al )\cdot  ( T_\chs ( \al \circ (-1)))=$$

$$ -   ( \al (\si )\circ (-1) )\cdot    E_{1,\chs}
-( \al (\si ) \circ (-1)) \cdot \frac{ (1- \chs )}{2}\de _0 +
$$

$$
 T_\chs ( \al )\cdot \Big(    E_{1,\chs} +  \Mc (\chs )+(1-\chs )\de _0+\frac{1- \chs }{2}\de _\chs \Big)+
$$

$$
 (-1)\Big(    E_{1,\chs} +  \Mc (\chs )+ (1-\chs )\de _0       +\frac{1- \chs }{2}  \de _\chs \Big)\cdot    E_{1,\chs} +
$$

$$
 (-1)\Big(    E_{1,\chs} +  \Mc (\chs )+ (1-\chs )\de _0       +\frac{1- \chs }{2}  \de _\chs \Big)\cdot   \frac{1- \chs }{2}  \de _0\, .
$$

\noindent
after applying several times Lemma 8.3.
We shall gather the remaining terms, which clearly form measures. These are terms

$$
 {\frac{(1-\chi (\si ) )^2}{8}}Y_0 ^2 +\sbp \ga _{ -b }^{(n)} (\si ) \big(  \sap Y_a\big) Y_b-  \sap \ga _{ -a }^{(n)} (\si )Y_a \big(  \sbp Y_b\big)  +
$$

$$
  (-1)    \sum _{b=1}^{p^n-1} \al _{p^n-b} ^{(n)} (\si )  \big( \sum _{b<a\leq p^n}Y_a\big) Y_b+ \sum _{a=1}^{p^n-1} \al _{p^n-a} ^{(n)} (\si )Y_a \big( \sum _{a<b\leq p^n}Y_b\big)  +
$$

$$
  {\frac{1}{2}}  \big(   {\frac{1-\chi (\si )}{2}} \big) ^2(Y_{\ls })^2+
$$

$${\frac{1-\chi (\si )}{2}}Y_0 \cdot         ( \sum _{b=0}^{p^n-1} \al _b ^{(n)} (\si )Y_ b)+ ({\frac{1-\chi (\si )}{2}}Y_{\ls }  )\cdot              (\sum _{b=0}^{p^n-1} \al _b ^{(n)} (\si )Y_ b) +
$$

$$
(-1) (  \sum _{a=0}^{p^n-1} \al _{p^n-a} ^{(n)} (\si )Y_ a )\cdot    {\frac{1-\chi (\si )}{2}}Y_0  +(  \sum _{a=0}^{p^n -1}\al _{\ls -a}^{(n)}(\si ) Y_a     )\cdot    {\frac{1-\chi (\si )}{2}}Y_0    +
$$

$$
 ({\frac{1-\chi (\si )}{2}}Y_{\ls }  )\cdot             {\frac{1-\chi (\si )}{2}}Y_0  -
  ( \sap \al _{a- \ls }^{(n)}(\si )  Y_a )\cdot      {\frac{1-\chi (\si )}{2}}Y_0    +
$$

$$
(-1)  ({\frac{1-\chi (\si )}{2}}Y_{\ls }  )\cdot          (  \sum _{b=0}^{p^n-1} \al _{p^n-b} ^{(n)} (\si )Y_ b  )      + ({\frac{1-\chi (\si )}{2}}Y_{\ls }  )\cdot        (  \sum _{b=0}^{p^n -1}\al _{\ls -b}^{(n)}(\si ) Y_b     )   +
$$

$$
( -1) ( \sap \al _{a- \ls }^{(n)}(\si )  Y_a )\cdot        ({\frac{1-\chi (\si )}{2}}Y_{\ls }  )    \, .
$$

From these terms we get the following sum of measures

\smallskip

\begin{equation}
\frac{  (1- \chs )^2  }{8} \de _0\cdot \de _0 +D_2(\si ) +\frac{  (1- \chs )^2 }{8} \de _\chs \cdot \de _\chs +\frac{  (1- \chs ) }{2} \de _0\cdot \al (\si )+
\tag{8.4.3}
\end{equation}

\smallskip

 \noindent
$$\frac{  (1- \chs ) }{2} \de _\chs \cdot \al (\si ) - \big( \al (\si )\circ (-1)\big)\cdot \frac{  (1- \chs ) }{2} \de _0 + \big( T_\chs ( \al (\si )\circ (-1))\big) \cdot \frac{  (1- \chs ) }{2} \de _0 +
$$

\smallskip

\noindent
$$\frac{  (1- \chs ) }{2} \de _\chs \cdot \frac{  (1- \chs ) }{2} \de _0 - ( T_\chs ( \al (\si )) \cdot \frac{  (1- \chs ) }{2} \de _0 -
\frac{  (1- \chs ) }{2} \de _\chs \cdot  (\al (\si )\circ (-1))+
$$

\smallskip

\noindent
$$\frac{  (1- \chs ) }{2} \de _\chs \cdot  \big( T_\chs ( \al (\si )\circ (-1))\big) - ( T_\chs ( \al (\si )) \cdot \frac{  (1- \chs ) }{2} \de _\chs =$$

\smallskip

$$
\frac{(-3)  (1- \chs )^2  }{8}    \de _0\cdot \de _0 +D_2(\si ) +\frac{  (1- \chs )^2 }{8} \de _\chs \cdot \de _\chs +\frac{  1- \chs  }{2} \de _0\cdot \al (\si ) + 
$$
$$
(-1) \big( \al (\si )\circ (-1)\big)\cdot \frac{  1- \chs  }{2} \de _0 +\frac{  (1- \chs )^2 }{4} \de _\chs \cdot \de _0 +\frac{  1- \chs  }{2} \de _\chs \cdot   T_\chs ( \al (\si )\circ (-1))+
$$
$$
(-1) T_\chs ( \al (\si )) \cdot \frac{  1- \chs  }{2} \de _\chs +   \frac{  1- \chs  }{2} \de _\chs\cdot   E_{1,\chs} -  E_{1,\chs}\cdot  \frac{  1- \chs  }{2} \de _0 +
$$
$$
(-1)\Mc (\chs )\cdot  \frac{  1- \chs  }{2} \de _0\, .
$$

Now we shall write terms which give products with the measure $\Mc (\chs )$. Then we have the following terms

$$
 {\frac{1-\chi (\si )}{2}} Y_{\ls }( \sum _{i=1}^{\ls - 1}Y_i)+  ({\frac{1-\chi (\si )}{2}}Y_{\ls }  )\cdot        {\frac{\chs - \ls}{ p^n}}(\sbp Y_b)   +
$$

$$
(-1)  ( \sap \al _{a- \ls }^{(n)}(\si ) Y_a) ( \sum _{b=1}^{\ls -1}Y_b ) - ( \sap \al _{a- \ls }^{(n)}(\si )  Y_a )\cdot        {\frac{\chs - \ls}{ p^n}}(\sbp Y_b)             +
$$

$$ 
 (  {\frac{\chs - \ls}{ p^n}}(\sap Y_a) )\cdot  ( \sum _{b=0}^{p^n-1} \al _b ^{(n)} (\si )Y_ b)+ (\sum _{a=1}^{\ls -1}Y_a)\cdot       ( \sum _{b=0}^{p^n-1} \al _b^{(n)} (\si )Y_ b)+
$$

$$
 (  {\frac{\chs - \ls}{ p^n}}(\sap Y_a) )\cdot        {\frac{1-\chi (\si )}{2}}Y_0     +
 (\sum _{a=1}^{\ls -1}Y_a)\cdot      {\frac{1-\chi (\si )}{2}}Y_0 +
 $$
$$
(-1)  (  {\frac{\chs - \ls}{ p^n}}(\sap Y_a) )\cdot       (  \sum _{b=0}^{p^n-1} \al _{p^n-b} ^{(n)} (\si )Y_ b  )        
 -  (\sum _{a=1}^{\ls -1}Y_a)\cdot      (  \sum _{b=0}^{p^n-1} \al _{p^n-b} ^{(n)} (\si )Y_ b  )\, ,
 $$
which give the following sums of measures

\smallskip

\begin{equation}
  {\frac{1-\chi (\si )}{2}}\de _\chs \cdot \Mc (\chs )- T_\chs (\al (\si ))\cdot  \Mc (\chs )+ \Mc (\chs ) \cdot \al (\si) + 
\tag{8.4.4.a}
\end{equation}
$$
\Mc (\chs ) \cdot  {\frac{1-\chi (\si )}{2}}\de _0-  \Mc (\chs ) \cdot \al (\si)\circ (-1)\, .
$$

\smallskip

There are some more terms which could  give products with $\Mc (\si )$.  
These terms are

$$
 {\frac{1}{2}}\big(  {\frac{\chs - \ls}{  p^n}}\big) ^2 \big( \sa Y_aY_b \big)+ {\frac{1}{2}} \Big(  \sum _{i=1}^{\ls - 1}Y_i  \Big) ^2 +   {\frac{1}{2}}         (\sum _{a=1}^{\ls -1}Y_a)\cdot       {\frac{\chs - \ls}{ p^n}}(\sbp Y_b)   +
$$
$$  
(  \sum _{a=0}^{p^n -1}\al _{\ls -a}^{(n)}(\si ) Y_a     )\cdot        {\frac{\chs - \ls}{ p^n}}(\sbp Y_b)       \, .
$$
 However we must add some additional terms.     We add the terms  
$$
{\frac{1}{2}}\big(  {\frac{\chs - \ls}{  p^n}}\sum _{a=0}^{p^n -1}Y_a \big) \big( \sum _{b=1}^{\ls -1}Y_b\big) +\big(  \sum _{a=0}^{p^n -1}\al _{\ls -a}^{(n)}(\si ) Y_a    \big) \cdot  \big( \sum _{b=1}^{\ls -1}Y_b\big)\, .
$$

Then we have the following sum of measures
\smallskip

\begin{equation}
 \frac{1}{2} \Mc (\chs )\cdot \Mc (\chs ) + T_\chs \big( \al \circ (-1)\big) \cdot \Mc (\chs )\, .
\tag{8.4.4.b}
\end{equation}

Adding (8.4.4.a) and (8.4.4.b) we get

\begin{equation}
  {\frac{1-\chi (\si )}{2}}\de _\chs \cdot \Mc (\chs )- T_\chs (\al (\si ))\cdot  \Mc (\chs )+ \Mc (\chs ) \cdot \al (\si) + 
\tag{8.4.4}
\end{equation}

$$
\Mc (\chs ) \cdot  {\frac{1-\chi (\si )}{2}}\de _0-  \Mc (\chs ) \cdot (\al (\si)\circ (-1))\, +
$$
$$
 \frac{1}{2} \Mc (\chs )\cdot \Mc (\chs ) + T_\chs \big( \al \circ (-1)\big) \cdot \Mc (\chs )\, =
$$

$$
  - \frac{1}{2} \Mc (\chs )\cdot \Mc (\chs ) +\Mc (\si ) \cdot E_{1,\chs } -  E_{1,\chs }\cdot  \Mc (\chs ) +
$$
$$
 \Mc (\chs )\cdot (1 - \chs )\de _0 -  (1 - \chs )\de _0 \cdot  \Mc (\chs )\, 
$$
by Lemma 8.3.
\smallskip

Now we look for terms which should give the measure $ ( -1) T_{\chs ,\chs }(D_2(\si ))$. These terms, including the substraction of one term we have added to get the measure $   T_\chs \big( \al \circ (-1)\big) \cdot \Mc (\chs ) $,  are

$$
\sum _{a=1}^{\ls}\al _{\ls -a}^{(n)}(\si )   Y_a \big( \sum _{b=1}^{a-1}Y_b\big) -   \sum _{b=1}^{\ls}\al _{\ls -b}^{(n)}(\si )            \big( \sum _{a=1}^{b-1}Y_a\big) Y_b +  
$$

\[
(-1) \sum _{ a=\ls +   1}^{p^n-1}\al _{\ls -a}^{(n)}(\si )  Y_a \big( \sum _{b=a+1}^{p^n}Y_b \big)+ 
 \sum _{ b=\ls +   1}^{p^n-1}\al _{\ls -b}^{(n)}(\si ) \big(\sum _{a=b+1}^{p^n}Y_a\big) Y_b\Big) + 
\]
\[
(-1) \sum _{b=0}^{p^n-1} \ga _{\ls -b}^{(n)}(\si ) \big(\sum _{a=0}^{p^n-1}Y_a\big)Y_b + \sum _{a=0}^{p^n-1} \ga _{\ls -a}^{(n)}(\si )  Y_a  \big(\sum _{b=0}^{p^n-1}Y_b\big) \, +
\]
$$
 (\sum _{a=1}^{\ls -1}Y_a)      (  \sum _{b=0}^{p^n -1}\al _{\ls -b}^{(n)}(\si ) Y_b     ) - \big(  \sum _{a=0}^{p^n -1}\al _{\ls -a}^{(n)}(\si ) Y_a    \big)   \big( \sum _{b=1}^{\ls -1}Y_b\big) \, .
 $$

\smallskip

These terms give a measure
\begin{equation}
  ( -1) T_{(\chs ,\chs )}(D_2(\si )) +T_\chs (\al \circ (-1))\cdot \de _\chs -  \de _\chs \cdot T_\chs (\al \circ (-1)) \,.
\tag{8.4.5}
\end{equation}

\smallskip

The remaining terms are 
\[
  {\frac{\chs - \ls}{2  p^n}}  \big( \sum _{1\leq b<a\leq p^n}Y_aY_b \big) +
  ( -1)  {\frac{\chs - \ls}{2  p^n}} \big( \sum _{1\leq a<b\leq p^n}Y_aY_b \big) + 
\]

$$
(-1) {\frac{1-\chi (\si )}{2}}   ( \sum _{i=1}^{\ls - 1}Y_i) Y_{\ls }\,+ {\frac{1}{2}} \sum _{1\leq j < i< \ls }\big(  Y_iY_j - Y_jY_i \big)+ 
$$

$$
      {\frac{1}{2}}         (\sum _{a=1}^{\ls -1}Y_a)\cdot       {\frac{\chs - \ls}{ p^n}}(\sbp Y_b)   + (\sum _{a=1}^{\ls -1}Y_a)\cdot     ({\frac{1-\chi (\si )}{2}}Y_{\ls }  ) +
 $$

$$
(-1){\frac{1}{2}}\big(  {\frac{\chs - \ls}{  p^n}}\sum _{a=0}^{p^n -1}Y_a \big) \big( \sum _{b=1}^{\ls -1}Y_b\big) \, .
$$
One checks that one gets the measure 

\begin{equation}
 {\frac{1}{2}}  \Nc (\chs ) \, .   
\tag{8.4.6}
\end{equation}

\noindent
The summation of terms from (8.4.1) to  (8.4.6) gives the formula of the theorem.  \hpb

\medskip

\noindent 
{\bf Proposition 8.5. }{\it The relation 
$$
\al (\si )- \al (\si )\circ (-1)=E_{1,\chs }+{\frac{1-\chs }{2}}\de _0
$$
is the consequence of $\Zbb /(2)$ and $\Zbb /(3)$ symmetries of $\Pbb ^1_{\overline \Qbb }\setminus \{ 0,1,\infty \}$.}

\medskip

\noindent
{\bf Proof.} It is a consequence of 
 $\Zbb /(2)$ and $\Zbb /(3)$ symmetries of $\Pbb ^1_{\overline \Qbb }\setminus \{ 0,1,\infty \}$ that the coefficient of the power series
 $\log \big( E_0(\ffk  _{\pi _0}(\si ))(X,Y)\big)$ at $YX^{2i-1}$ is equal ${\frac{-B_{2i}}{2\cdot (2i)!}}(\chs ^{2i} -1)$ (see \cite[the proof of Prop. 3.1 and Prop. 3.1]{W14}).
Hence the coefficient of $ E_0(\ffk  _{\pi _0}(\si ))(X,Y)$ at $YX^{2i-1}$ is the same, as the power series  $ E_0(\ffk  _{\pi _0}(\si ))(X,Y)$ 
has no terms of degree $1$.
Hence it follows from \cite[Theorem 2.6]{W11} or \cite{NW} that
$$
{\frac{1}{(2i-1)!}}\int _{\Zbb _p} x^{2i-1}d\al (\si )(x)={\frac{-B_{2i}}{2\cdot (2i)!}}(\chs ^{2i} -1)\, .
$$
Therefore we have 
$$
\int _{\Zbb _p} x^{2i-1}d\al (\si )(x) ={\frac{-B_{2i}}{2\cdot 2i}}(\chs ^{2i} -1)\, .
$$
Observe that $\int _{\Zbb _p} x^kd\big( \al (\si )\circ (-1)\big)(x) =(-1)^k\int _{\Zbb _p} x^{k}d\al (\si )(x) $ for $k\geq 0$.
Hence it follows that 
$$
\int _ \zp x^{k-1}d\big( \al (\si ) - \al (\si )\circ (-1)\big) (x)={\frac{B_k}{k}}(1-\chs ^k)
$$
for $k>1$ and $\int _\zp d\big( \al (\si ) - \al (\si )\circ (-1)\big) (x)=0$. On the other side
$$
\int _\zp x^{k-1}d\big( E_{1,\chs }+{\frac{1-\chs }{2}}\de _0\big) (x)={\frac{B_k}{k}}(1-\chs ^k)
$$
(see \cite[Chapter 2, Theorem 2.3]{L}). Hence it follows that 
$$
F(\al (\si ) - \al (\si )\circ (-1))(X)=F (   E_{1,\chs }+{\frac{1-\chs }{2}}\de _0      )(X)\, ,
$$
and in a consequence $\al (\si ) - \al (\si )\circ (-1)=E_{1,\chs }+{\frac{1-\chs }{2}}\de _0$.
\hpb

\medskip

\noindent 
{\bf Proposition 8.6. }{\it  The relation  $\al (\si ) - \al (\si )\circ (-1)=E_{1,\chs }+{\frac{1-\chs }{2}}\de _0$ implies the relation 
$\al (\si ) -\al (\si )\circ (-1)+ T_{\chs } (\al (\si )\circ (-1) )    -   T_{\chs }(\al (\si )) +
 \Mc (\chs ) +{\frac{1-\chs }{2}}\de _0 +{\frac{1-\chs }{2}}\de _{\chs }=0$.}

\medskip

\noindent
{\bf Proof.} It follows from Lemma 8.3. ii)  that $T_{\chs } \big( E_{1,\chs }+{\frac{1-\chs }{2}}\de _0 \big) = E_{1,\chs } +  \Mc (\chs ) +
(1-\chs )\de _0 +{\frac{1-\chs }{2}}\de _{\chs }$.
Hence it follows that $\al (\si ) - \al (\si )\circ (-1) + T_{\chs } (\al (\si )\circ (-1) )    -   T_{\chs }(\al (\si )) =E_{1,\chs }+{\frac{1-\chs }{2}}\de _0 -
 E_{1,\chs } -  \Mc (\chs ) -
(1-\chs )\de _0 - {\frac{1-\chs }{2}}\de _{\chs }= -  \Mc (\chs ) -
{\frac{(1-\chs )}{2}}\de _0 -  {\frac{1-\chs }{2}}\de _{\chs } $. \hpb

\medskip

\noindent
 Remark.  The image of the map $\ffk _\pi :G_{\Qbb}\to  \pi _1(V,\01 )$ satisfies $\Zbb /(2)$, $\Zbb /(3)$ and $\Zbb /(5)$ relations.
The image of the composition ${\bf K} \circ \ffk _\pi :G_{\Qbb}\to      \prod _{m=0}^\infty   {\rm Meas} \big( (\zp )^m,{\frac{1}{m!}}\zp \big)$
satisfies the relations obtained from the octagonal relations on $\pi _1(V_n,\01 )$ for all $n$. It follows from Proposition 8.6 that there are some more relations.

\medskip

Let us set $\Lc _{-1}:=G_\Qbb$, $\Lc _0:= G_{\Qbb (\mu _{p^\infty })}$ and 
$\Lc _m:=\{ \si \in  G_{ \Qbb (\mu _{p^\infty })  }\mid  { \Kc}  _i(\si )=0$  $ {\rm for }\; 1\leq i \leq m \}$.
The filtration $\{ \Lc _m\} _{m=-1}^\infty$ of $G_\Qbb$ coincides with the filtration defined in \cite[page 97]{W10}
and denoted by the same symbols. The filtration $\{ \Lc _m\} _{m=-1}^\infty$ is the filtration by closed, normal subgroups of $G_\Qbb$. We left
the reader the proof of this fact (see also   \cite[Lemma 3.1.]{W10}).

\medskip

It would be very laborious to write explicitely the similar relations as in Theorem 8.4 for measures $\Kc _m(\si )$ for $m\geq 3$ and $\si \in G_\Qbb$. We have only the following result.
 
\medskip

\noindent 
{\bf Proposition 8.7. } (see also \cite[Theorem 4.2.]{W10} for $\si \in \Lc _{m-1}$) {\it  Let $\si \in G_\Qbb$ and let $m\geq   1$.
Let us set ${\overline \chs}:=(\chs, \chs ,\ldots ,\chs )\in (\zp )^m$.
There is a measure $h_m(\si )\in {\rm Meas} \big( (\zp )^m,{\frac{1}{m!}}\zp \big)$ such that
\begin{enumerate}
\item[i)] $h_m(\si )=0$  for $\si \in \Lc _{m-1}$, 
\item[ii)]  $\Kc _m (\si ) -\Kc _m(\si )\circ (-1) + T_{ {\overline \chs}}(\Kc _m(\si )\circ (-1))- T_{ {\overline \chs}}(\Kc _m(\si )) +h_m(\si )=0$ for
$\si \in G_{\Qbb }$.
\end{enumerate} }

\medskip

\noindent
{\bf Proof.} It follows from Theorem 2.7 and Lemmas 2.1, 2.2 and 2.3 for $N=p^n$ that comparing coefficients at $\bf Y_a$ (${\bf a}=(a_1,\ldots ,a_m)\in (\Zbb /(p^n))^m$, ${\bf Y_a}=Y_{a_1,n}Y_{a_2,n}\ldots Y_{a_m,n}$) after the embedding $E_n:\pi _1(V_n,\01 )\to \Qbb _p \{\{\Yc _n\}\}$
we get $\lam _{\bf Y_ a}^{(n)}(\si ) - \lam _{{\bf Y}_{{\bf -  a}}}^{(n)}(\si ) + \lam _{\bf Y _{{\overline \chs} -  a}}^{(n)}(\si )-  \lam _{  \bf   Y_{ a-{\overline \chs} }}^{(n)}(\si ) $
plus an other term, which we denote by $h_{\bf a}^{(n)}(\si )$.
It is clear that $
\lam _{\bf Y_ a}^{(n)}(\si ) - \lam _{{\bf Y}_{{\bf -  a}}}^{(n)}(\si ) + \lam _{\bf Y _{{\overline \chs} -  a}}^{(n)}(\si )-  \lam _{  \bf   Y_{ a-{\overline \chs} }}^{(n)}(\si ) 
+h_{\bf a}^{(n)}(\si )=0$
as on the left hand side of the octagonal relation we have $1$. Moreover it is also clear that $\big( \big( \Zbb /(p^n)\big)^m \ni {\bf a}\mapsto h_{\bf a}^{(n)}(\si )\big) _{n\in \Nbb }$ is a measure, as $\Kc _m(\si )$ is a measure and the sum is $0$, hence zero measure on $(\zp )^m$.

The term $h_{\bf a}^{(n)}(\si )$ is formed with coefficients which vanish for $\si \in G_{\Qbb (\mu _{p^\infty})}$ as for example $\frac{\langle \si \rangle _n -\chs }{p^n}$ and  with coefficients which are integrals against measures $\Kc _i (\si )$ for $i<m$ (coefficients at monomials of degree $m$ containing $X_n$). Hence it follows that $h_{\bf a}^{(n)}(\si )=0$ for $\si \in \Lc _{m-1}$. \hpb

\medskip

\noindent
Remark.  One should be able to show that $\be _m$ is non-zero on $\Lc _{m-1}$.

\medskip

\noindent 
{\bf Proposition 8.8.} {\it We have
$$
\bigcap _{m=-1}^\infty  \Lc _m=\ker (G_\Qbb \to \Aut (\pi _1(V,\01)).
$$
}

\noindent
{\bf Proof.} The corollary follows from Proposition 7.2.    \hpb

\bigskip

\section{Towards $p$-adic multiple  zeta functions } 
\smallskip

In \cite{W11} we asked how to get from the measures $\Kc  _n(\si )$ a $p$-adic multiple zeta functions. We propose here a construction which we hope is a step in this direction.

 Let $S_m$  be the group of permutations of the set $\{ 1,2,\ldots ,m\}$.  We shall define actions of the groups $S_m$ and $\{1,-1\}^m$ on the $\zp$-module $\zp [[(\zp )^m]]$. Let $s\in S_m$. Then we view $s$ as a permutation matrix. Hence $s$ induces a linear isomorphism $s:(\zp )^m \to (\zp )^m$,
which we extend by linearity and continuity to an isomorphism $s: \zp [[(\zp )^m]] \to \zp [[(\zp )^m]]$. Observe also that $s$ induces an isomorphism $s:(\Zbb /(p^n))^m \to (\Zbb /(p^n))^m $ for any $n$.

\medskip

\noindent 
{\bf Proposition 9.1. }{\it   Let $\si \in G_\Qbb$ and let $m\geq 2$.   Let  ${\bf a}=(a_1,a_2, \ldots ,a_m)\in \big( \Zbb /(p^n) \big)^m$. Then we have:
\begin{enumerate}
\item[i)]
\[
\sum _{s\in S_m}\Kc  _m ^{(n)} (\si )   (s({\bf a}) )  =\Kc  _1 ^{(n)} (\si )  (a_1) \cdot  \Kc  _1 ^{(n)} (\si )  (a_2)  \cdot \ldots \cdot  \Kc  _1 ^{(n)} (\si )  (a_m)  \, ,
\]
\item[ii)] 
\[
\sum _{s\in S_m} s( \Kc _m(\si ))=\Kc _1 (\si )\cdot \Kc _1 (\si )\cdot \ldots \cdot \Kc _1 (\si )=\Kc _1 (\si )^m.
\]
\end{enumerate}
}

\noindent
{\bf Proof.}  The point i) of the proposition follows from the definition of the embedding

\noindent
 $E_n:\pi _1(V_n,\01 )\to  \Qbb _p\{\{ \Yc _n \}\}$ and from the fact that the elements of the form $\exp w$, where $w$ is a Lie element, satisfy the shuffle relations (see \cite[Theorem 2.5]{R}). The point ii) follows from the point i).
\hpb

\medskip

\noindent
Let $\varepsilon =(\varepsilon _1, \varepsilon  _2, \ldots  ,   \varepsilon _m)\in (\{1,-1\})^m$  and let $(x_1,x_2,\ldots , x_m) \in (\zp )^m$.
We set 
\[
\varepsilon ((x_1,x_2,\ldots , x_m)):=(\prod _{i=1}^m\varepsilon _i  ) (\varepsilon _1  x_1, \varepsilon  _2 x_2, \ldots  ,   \varepsilon _m x_m)\,.
\]
We extend by linearity and by continuity to the action of the group  $(\{1,-1\})^m$ on the $\zp$-module $\zp [[(\zp)^m]]$. The actions of $S_m$ and  $(\{1,-1\})^m$ give together an action of the semi-direct product $(\{1,-1\})^m  \rtimes S_m$ on  the $\zp$-module $\zp [[(\zp)^m]]$.

\medskip

\noindent 
{\bf Theorem  9.2. }{\it    Let $\si \in G_\Qbb$ and let $m\geq 2$.  We have
\[
\sum _{g\in (\{1,-1\})^m  \rtimes S_m} g(\Kc  _m(\si ))=\big( E_{1,\chs }+ {\frac{1-\chs}{2}}\de  _0\big) ^m \,
\] 
in $\zp [[(\zp )^m]]$.}

\noindent
{\bf Proof.} Let us fix ${\bf a}=(a_1,a_2, \ldots ,a_m)\in \big(\Zbb /(p^n)\big) ^m$. Let  ${\bf a ^\prime}=(-a_1,a_2, \ldots ,a_m)$. Then we have
\[
\sum _{s\in S_m}\Kc _m ^{(n)}(\si )  (s({\bf a}) )  =\Kc _1 ^{(n)} (\si )  (a_1) \cdot  \Kc _1 ^{(n)} (\si )  (a_2)  \cdot \ldots \cdot  \Kc _1 ^{(n)} (\si )  (a_m)    \;{\rm and}\; 
\]
\[
\sum _{s\in S_m}\Kc _m ^{(n)}(\si )  (s({\bf a ^\prime}) )  =\Kc _1 ^{(n)} (\si )  (- a_1) \cdot  \Kc _1 ^{(n)} (\si )  (a_2)  \cdot \ldots \cdot  \Kc _1 ^{(n)} (\si )  (a_m)   \, .
\]
by Proposition 9.1. Hence it follows that
\[
\sum _{s\in S_m}        \Kc _m ^{(n)}(\si )  (s({\bf a}) )        -   \sum _{s\in S_m}        \Kc _m ^{(n)}(\si )  (s({\bf a ^\prime}) )            =
\]
\[
\big(  E_{1,  \chs }^{(n)}(\al _1)+ {\frac{1-\chs}{2}}\de  _0 (\al _1)\big)  \cdot       \Kc _1 ^{(n)} (\si )  (a_2)               \cdot \ldots \cdot   \Kc _1 ^{(n)} (\si )  (a_m) 
\]
by Lemma 8.3. Repeating the same argumets $m$-times we get the proposition.  \hpb

\medskip

\noindent 
{\bf Corollary  9.3. }{\it Let $m=2$ and let $\si \in G_\Qbb $.  Then  we have the following identity in $\zp [[T_1,T_2]]$
\[
P_{\Kc _2(\si )} (T_1,T_2)+P_{\Kc _2(\si )} (T_2,T_1)- P_{\Kc _2(\si )} ((1+T_1)^{-1},T_2) - P_{\Kc _2(\si )} ( T_2, (1+T_1)^{-1})  
\]
\[
 - P_{\Kc _2(\si )} ( T_1, (1+T_2)^{-1}) - P_{\Kc _2(\si )} ((1+T_2)^{-1},T_1) +
\]
\[
 P_{\Kc _2(\si )} ((1+T_1)^{-1},(1+T_2)^{-1})  + P_{\Kc _2(\si )} ((1+T_2)^{-1},(1+T_1)^{-1}) =
\]
\[
\big({\frac{1}{T_1}}-  {\frac{\chs}{(1+T_1)^{\chs } - 1}}+{\frac{1-\chs}{2}}\big) \cdot \big({\frac{1}{T_2}}-  {\frac{\chs}{(1+T_2)^{\chs } - 1}}+{\frac{1-\chs}{2}}\big) \,.
\] \hpb
}

\medskip

\noindent 
 Remark.   In \cite[Definition 2.9.]{F} the authors use the measures $\big(  E_{1,  \chs }\big) ^m$ to construct $p$-adic multiple $L$-functions so we hope Theorem 9.2 recovers correctly the measures which give $p$-adic multiple zeta functions.

\bigskip

\section{Correction of two results of \cite{W10} }
\smallskip
 
We shall correct two results \cite[Corollaries 5.2 and 5.3]{W10}.

\medskip

\noindent 
{\bf Lemma 10.1.}  (Correction of \cite[Corollary 5.2]{W10}) {\it Let $n$ and $r$ be positive integers. Let $n_0,n_1,\ldots ,n_r$ be non-negative integers. Let $0\leq i_k<p^n$ for $1\leq k\leq r$. 
Let $m=\sum _{i=0}^rn_i$. Let us set ${\bf i} =(i_1,\ldots ,i_r)$, ${\bf 1}=(1,\ldots ,1)\in (\Zbb )^r$ and ${\bf p^n}=p^n\bf 1$.
For ${\bf a}=(a_1,\ldots ,a_r)\in (\Zbb )^r$ we set 
$$
F_{\bf a}(x_1,\ldots ,x_r):=({\frac{x_1-x_2-a_1+a_2}{p^n}})^{n_1}\ldots 
({\frac{x_{r-1}-x_r-a_{r-1}+a_r}{p^n}})^{n_r} .
$$
 Let $\si \in \Lc _{r-1}$. Then we have
$$
\int _{{\bf i}+p^n (\zp )^r}(\frac{i_1-x_1}{p^n})^{n_0}F_{\bf i}(x_1,\ldots ,x_r)(\frac{x_r-i_r}{p^n})^{n_r}d\Kc _r (\si )(x_1,\ldots ,x_r) \, +
$$
$$
(-1)^{m+1} \int _{ {\bf p^n}- {\bf i}+p^n (\zp )^r}(\frac { (p^n -i_1)-x_1}{p^n}-1)^{n_0}F_{ {\bf p^n}-   \bf i}(x_1,\ldots ,x_r)(\frac{x_r-(p^n -i_r)}{p^n}+1)^{n_r} 
$$
$d\Kc _r (\si )(x_1,\ldots ,x_r) \, +$
$$
(-1)^{m} \int _{ {\bf p^n}+{\bf 1}- {\bf i}+p^n (\zp )^r}(\frac { (p^n -i_1+1)-x_1}{p^n}-1)^{n_0}F_{ {\bf p^n}+{\bf 1}-   \bf i}(x_1,\ldots ,x_r)(\frac{x_r-(p^n -i_r+1)}{p^n}+1)^{n_r} 
$$
$d\Kc _r (\si )(x_1,\ldots ,x_r) \, +$
$$
(-1)\int _{{\bf i}-{\bf 1}+p^n (\zp )^r}(\frac{i_1-1-x_1}{p^n})^{n_0}F_{{\bf i}-{\bf 1}}(x_1,\ldots ,x_r)(\frac{x_r-(i_r-1)}{p^n})^{n_r}d\Kc _r (\si )(x_1,\ldots ,x_r) =0\, .
$$}

\noindent
{\bf Proof.} Let $\mu =\sum _{i=1}^{r-1}n_i$ and let $\al \in \Zbb$. Observe that
$F_{\bf a}(-x_1,\ldots ,-x_r)=(-1)^\mu F_{-\bf a}(x_1,\ldots ,x_r)$,
$F_{\bf a}(1-x_1,\ldots ,1-x_r)=(-1)^\mu F_{{\bf 1}-\bf a}(x_1,\ldots ,x_r)$
and

\noindent
$F_{{\bf a}+\al \bf1}(x_1,\ldots ,x_r)=F_{\bf a}(x_1,\ldots ,x_r)$.

Observe also that for any ${\bf v}\in (\zp )^r$ we have ${\bf a}+p^n{\bf v}+p^n (\zp )^r ={\bf a}+p^n (\zp )^r$ and
${\bf a}-p^n (\zp )^r= {\bf a}+p^n (\zp )^r$.

Making a change of variables $y_i=-x_i$ ($1\leq i\leq r$) we get

\smallskip

\begin{equation}
\int _{{\bf i}+p^n (\zp )^r}(\frac{i_1-x_1}{p^n})^{n_0}F_{\bf i}(x_1,\ldots ,x_r)(\frac{x_r-i_r}{p^n})^{n_r}d\Kc _r (\si )(-x_1,\ldots ,-x_r) =
\tag{10.a}
\end{equation}
$$
\int _{-{\bf i}-p^n (\zp )^r}(\frac{i_1+y_1}{p^n})^{n_0}F_{\bf i}(-y_1,\ldots ,-y_r)(\frac{-y_r-i_r}{p^n})^{n_r}d\Kc _r (\si )(y_1,\ldots ,y_r) =
$$
$$
(-1)^m \int _{-{\bf i}+p^n (\zp )^r}(\frac{-i_1-y_1}{p^n})^{n_0}F_{-\bf i}(y_1,\ldots ,y_r)(\frac{y_r+i_r}{p^n})^{n_r}d\Kc _r (\si )(y_1,\ldots ,y_r) =
$$
$
(-1)^{m}\int _{{\bf p^n}-{\bf i}+p^n (\zp )^r}(\frac{(p^n -i_1)-x_1}{p^n}-1)^{n_0}F_{{\bf p^n}-\bf i}(x_1,\ldots ,x_r)(\frac{x_r-(p^n -i_r)}{p^n}+1)^{n_r}$

\noindent
$d\Kc _r (\si )(x_1,\ldots ,x_r) .$

\smallskip

Similarly after changes of variables $y_i=1-x_i$ ($1\leq i\leq r$) and $y_i=x_i-1$ ($1\leq i\leq r$) we get 
\smallskip

\begin{equation}
\int _{{\bf i}+p^n (\zp )^r}(\frac{i_1-x_1}{p^n})^{n_0}F_{\bf i}(x_1,\ldots ,x_r)(\frac{x_r-i_r}{p^n})^{n_r}d\Kc _r (\si )(1-x_1,\ldots ,1-x_r) =
\tag{10.b}
\end{equation}
$\noindent
(-1)^{m} \int _{ {\bf p^n}+{\bf 1}- {\bf i}+p^n (\zp )^r}(\frac { (p^n -i_1+1)-x_1}{p^n}-1)^{n_0}F_{ {\bf p^n}+{\bf 1}-   \bf i}(x_1,\ldots ,x_r)(\frac{x_r-(p^n -i_r+1)}{p^n}+1)^{n_r}$

$d\Kc _r (\si )(x_1,\ldots ,x_r) 
$

\smallskip

and

\smallskip

\begin{equation}
\int _{{\bf i}+p^n (\zp )^r}(\frac{i_1-x_1}{p^n})^{n_0}F_{\bf i}(x_1,\ldots ,x_r)(\frac{x_r-i_r}{p^n})^{n_r}d\Kc _r (\si )(x_1-1,\ldots ,x_r-1) =
\tag{10.c}
\end{equation}
$
\int _{{\bf i}-{\bf 1}+p^n (\zp )^r}(\frac{i_1-1-x_1}{p^n})^{n_0}F_{{\bf i}-{\bf 1}}(x_1,\ldots ,x_r)(\frac{x_r-(i_r-1)}{p^n})^{n_r}d\Kc _r (\si )(x_1,\ldots ,x_r) \, .
$

\smallskip

Let $\si \in \Lc _{r-1}$. Then it follows from Proposition 8.7  (see also \cite[Theorem 4.2.]{W10}) that
$d(\Kc _r (\si )(x_1,\ldots ,x_r) - \Kc _r (\si )(-x_1,\ldots ,-x_r) + \Kc _r (\si )(1-x_1,\ldots ,1-x_r) - \Kc _r (\si )(x_1-1,\ldots ,x_r-1) )=0$.
Hence the lemma follows from the equalities (10.a), (10.b) and (10.c).   \hpb

\medskip
Let $n$ and $r$ be positive integers. Let $n_0,n_1,\ldots ,n_r$ be non-negative integers.
Let ${\bf a}=(a_1,\ldots ,a_r)\in (\Zbb /(p^n) )^r$. We set
$$
w({\bf a})=Y_{a_1,n}X_n^{n_1}Y_{a_2,n}X_n^{n_2}\ldots Y_{a_{r-1},n}X_n^{n_{r-1}}Y_{a_r,n}\,.
$$
We recall that $G_\Qbb$ has a (weight) filtration associated with the action of $G_\Qbb$ on $\pi _1(V_n,\01 )$. The filtration is defined in the following way:
$$
H^{(n)}_0:=G_\Qbb , H^{(n)}_1:=G_{\Qbb (\mu _{p^\infty })} , H^{(n)}_k:= \{\si \in G_{\Qbb (\mu _{p^\infty })} \mid E_n(\ffk _{\pi _n}(\si ))\equiv 1
\;{\rm mod}\; I_n^k\}
$$
for $k>1$.

\medskip

\noindent 
{\bf Proposition 10.2. }(Correction of \cite[Corollary 5.3]{W10}) {\it The assumptions and notations are as in Lemma 10.1.
\begin{enumerate}
\item[i)] Let $\si \in \Lc _{r-1}$. If $n_0=n_r=0$ and $1<i_k<p^n$ for $1\leq k \leq r$ then
$$
\lam ^{(n)}_{w({\bf i})}(\si )- (-1)^m\lam ^{(n)}_{w({-\bf i})}(\si )+(-1)^m \lam ^{(n)}_{w({{\bf 1}-\bf i})}(\si )-\lam ^{(n)}_{w({\bf i}-{\bf 1})}(\si )=0.
$$
\item[ii)] Let $\si \in \Lc _{r-1}\cap H^{(n)}_{m+r-1}$. Then
$$
\lam ^{(n)}_{X_n^{n_0}w({\bf i})X_n^{n_r}   }(\si )- (-1)^m\lam ^{(n)}_{ X_n^{n_0} w({-\bf i}) X_n^{n_r} }(\si )+(-1)^m \lam ^{(n)}_{ X_n^{n_0} w({{\bf 1}-\bf i})  X_n^{n_r} }(\si ) +
$$
$(-1)\lam ^{(n)}_{  X_n^{n_0} w({\bf i}-{\bf 1})  X_n^{n_r}  }(\si )=0.$
\end{enumerate}
}

\smallskip

\noindent
{\bf Proof.}  First we show the point i). Observe that if $1<i_k<p^n$ for $1\leq k \leq r$ then $0\leq p^n- i_k <p^n$, $0\leq p^n +1-i_k<p^n$ and
$0\leq i_k -1<p^n$ for $1\leq k \leq r$. Hence it follows from Corollary 10.1 and Theorem 3.1 that $\lam ^{(n)}_{w({\bf i})}(\si )- (-1)^m\lam ^{(n)}_{w({-\bf i})}(\si )+(-1)^m \lam ^{(n)}_{w({{\bf 1}-\bf i})}(\si )-\lam ^{(n)}_{w({\bf i}-{\bf 1})}(\si )=0$.

In order to show the point ii) we shall also use the formula of Lemma 10.1. We need to consider more carefully the last three integrals of the formula.
We have 
\begin{equation}
 \int _{ {\bf p^n}- {\bf i}+p^n (\zp )^r}(\frac { (p^n -i_1)-x_1}{p^n}-1)^{n_0}F_{ {\bf p^n}-   \bf i}(x_1,\ldots ,x_r)(\frac{x_r-(p^n -i_r)}{p^n}+1)^{n_r} 
\tag{10.d}
\end{equation}

\noindent
$d\Kc _r (\si )(x_1,\ldots ,x_r) =$

\noindent
$
 \int _{ {\bf p^n}- {\bf i}+p^n (\zp )^r}(\frac { (p^n -i_1)-x_1}{p^n})^{n_0 }F_{ {\bf p^n}-   \bf i}(x_1,\ldots ,x_r)(\frac{x_r-(p^n -i_r)}{p^n})^{n_r}d\Kc _r (\si )(x_1,\ldots ,x_r) +
$

\noindent
$
\sum _{\al =0}^{n_0}\sum _{\be =0, (\al ,\be )\neq (0,0)}^{n_r}  (-1)^\al \int _{ {\bf p^n}- {\bf i}+p^n (\zp )^r}(\frac { (p^n -i_1)-x_1}{p^n})^{n_0 -\al}F_{ {\bf p^n}-   \bf i}(x_1,\ldots ,x_r)
$

\noindent
$
(\frac{x_r-(p^n -i_r)}{p^n})^{n_r-\be}d\Kc  _r (\si )(x_1,\ldots ,x_r)
$.

\medskip

Let us assume that $0<i_k<p^n$ for $1\leq k \leq r$. Then we have $0\leq p^n - i_k<p^n$ for $1\leq k \leq r$. It follows from Theorem 3.1 that the integral (10.d) is equal 
$(\prod _{i=0}^rn_i! ) \lam ^{(n)}_{ X_n^{n_0} w({-\bf i}) X_n^{n_r} }(\si )$ $+$ sums of rational multiples of $\lam ^{(n)}_{ X_n^{\al } w({-\bf i}) X_n^{\be} }(\si )$
with $\al +\be<n_0+n_r$. Observe that for $\si \in H^{(n)} _{m+r-1}$ the terms $\lam ^{(n)}_{ X_n^{\al } w({-\bf i}) X_n^{\be} }(\si )$ with $\al +\be<n_0+n_r$ vanish.
All other cases are dealt similarly.
\hpb

\bigskip

\section{Examples in the Betti - de Rham context } 
\smallskip

In this section we briefly  indicate some comparisons with the complex (Betti - de Rham) situation. One can find some examples of such comparisons in 
\cite{NW2}, in particular \cite[pages 276 and 277]{NW2} and in \cite{NW3} and one easily notice some discrepancies.
 In the complex situation, the $\Cbb$-algebra $ \Cbb \{\{ X,Y\}\}$ is on generators $X=\big( \frac{dz}{z} \big)^*$ and  $Y=\big( \frac{dz}{z-1} \big)^*$  - the dual base of the base of one-forms with logarithmic singularities on $\Pbb ^1(\Cbb )\setminus \{0,1,\infty \}$. Then one have $k_*(X)=-X-Y$ and $k_*(Y)=Y$ for $k(\zfk )=1/\zfk$. However this does not holds in the $\Qbb _p$-algebra $\qp \{ \{ X,Y\} \}$ with the embedding considered in this paper.
The same remarks concerns also $\Pbb ^1(\Cbb )\setminus \big(\{0,\infty \}\cup \mu _{p^n}\big)$ and morphisms of it. Perhaps one can improve the situation on the Galois side if one chooses suitables embeddings of fundamental groups into $\qp$-algebras. We finish this section with two examples.

\medskip

\noindent 
{\bf Example  11.1. }  The elements
$$
\sum _{k=0}^{N-1}\big( \int _\01 ^{ \frac{1}{N} \10 } {\frac{dz}{z-\xi_N^k}}\big) [k]\in \Cbb [\Zbb/(N)]
$$
(the integration is along the path $\pi$ from section 1) for $N\in \Nbb$ form a distribution on $\hat \Zbb$  (see \cite[Theorem 4]{NW3}), which we denote by $\Zc$. Let $0<k<N$.
Then
$$
 \int _\01 ^{ \frac{1}{N} \10 } {\frac{dz}{z-\xi_N^k}}= \int _0 ^1{\frac{dz}{z-\xi_N^k}}= \int _{-\xi _N^k} ^{ 1-\xi _N^k  }{\frac{dz}{z }}={\frac{1}{2}}\pi i-\frac{k}{N} \pi i\, .
$$
Hence for any $0<k<N$,
$$
 \int _\01 ^{ \frac{1}{N} \10 } {\frac{dz}{z-\xi_N^k}}- \int _\01 ^{ \frac{1}{N} \10 } {\frac{dz}{z-\xi_N^{-k}}}=-2\pi i \big( \frac{k}{N} - \frac{1}{2}\big) =
-2\pi iE_1^{(N)}\big( \frac{k}{N}\big)\, .
$$
Hence it follows that
$$
\Zc - \Zc \circ (-1)=-2\pi i E_1 -\pi i \de _0\, .
$$
One can compare this result with the identity 
$$
\Kc _1(\si ) - \Kc _1(\si ) \circ (-1)=E_{1,\chi (\si )}+\frac{1-\chi (\si)}{2}\de _0
$$
(see \cite[Proposition 7]{W12}).

The octagonal relations in the complex (Betti - de Rham) case imply the following relation between distributions on $\hat \Zbb$,
$$
\Zc - \pi i \de _0 - \Zc \circ (-1) -2\pi i \Hc + T_1 ( \Zc \circ (-1) )-\pi i \de _1 - T_1(\Zc ) =0\, ,
$$
if one itegrates along the paths on Picture 2  the one-form $\frac{dz}{z }\otimes (\frac{dz}{z })^* +\sum _{k=0}^{N-1}\frac{dz}{z-\xi_N^k}\otimes (\frac{dz}{z-\xi_N^k})^*$ and where $\Hc$ is a distribution on $\hat \Zbb$ given by $\Hc ^{(N)}:\Zbb /(N)\to \Cbb,\; [k]\mapsto \frac{1}{N} $ for $N\in \Nbb$. The above formula is the analogue of the formulas of Proposition 8.1.

\medskip

\noindent 
{\bf Example  11.2. }  Let $0<k<N$. Then we have
$$
\int _\01 ^{ \frac{1}{N} \10 } {\frac{dz}{z-\xi_N^k}},\frac{dz}{z } + \int _\01 ^{ \frac{1}{N} \10 } {\frac{dz}{z-\xi_N^{-k}}},\frac{dz}{z } =\frac{1}{2}(2\pi i)^2 B_2(\frac{k}{N})\, .
$$
We left the stated formula without proof. However one can see also \cite[(6.26)]{NW2}. The above formula should be compared with the formula of Corollary 3.3.

\vglue 2cm

\vglue 1cm

\vglue 1cm

\noindent 179 Piste de L'Uesti

\noindent  06910 Pierrefeu,

\noindent France

\smallskip

\noindent {\it E-mail address} zdzislaw.wojtkowiak@gmail.com

\medskip

\end{document}